\documentclass[leqno,11pt,a4paper]{amsart}
\pdfoutput=1

\usepackage[full]{textcomp}
\usepackage{newpxtext}
\usepackage{diagbox}
\usepackage{comment}
\usepackage{cabin}
\usepackage[varqu,varl]{inconsolata}
\usepackage[bigdelims,vvarbb]{newpxmath}
\usepackage[cal=cm,bb=esstix,bbscaled=1.05]{mathalfa}
\usepackage[margin=2.395cm]{geometry}
\usepackage{savesym}
\let\savedegree\bigtimes
\let\bigtimes\relax
\usepackage{mathabx}
\usepackage{amsmath}
\usepackage{graphicx}
\let\bigtimes\savedegree

\usepackage[dvipsnames]{xcolor}
\usepackage{tikz}
\usepackage{tikz-cd}
\usepackage{cite}
\usepackage{enumitem}
\usepackage{amsfonts}

\setlist[enumerate]{labelsep=*, leftmargin=1.5pc}
\setlist[enumerate]{label=\normalfont(\roman*), ref=\roman*}

\usepackage{aliascnt}

\newtheorem{thm}{Theorem}[section]

\newaliascnt{lemma}{thm}
\newtheorem{lemma}[lemma]{Lemma}
\aliascntresetthe{lemma}

\newaliascnt{cor}{thm}
\newtheorem{cor}[cor]{Corollary}
\aliascntresetthe{cor}

\newaliascnt{prop}{thm}
\newtheorem{prop}[prop]{Proposition}
\aliascntresetthe{prop}

\newtheorem{innercustomthm}{Theorem}
\newenvironment{customthm}[1]
  {\renewcommand\theinnercustomthm{#1}\innercustomthm}
  {\endinnercustomthm}

\theoremstyle{definition}

\newaliascnt{example}{thm}
\newtheorem{example}[example]{Example}
\aliascntresetthe{example}

\newaliascnt{notation}{thm}
\newtheorem{notation}[notation]{Notation}
\aliascntresetthe{notation}

\newaliascnt{conj}{thm}
\newtheorem{conj}[conj]{Conjecture}
\aliascntresetthe{conj}

\newaliascnt{question}{thm}
\newtheorem{question}[question]{Question}
\aliascntresetthe{question}

\newaliascnt{remark}{thm}
\newtheorem{remark}[remark]{Remark}
\aliascntresetthe{remark}

\newaliascnt{definition}{thm}
\newtheorem{definition}[definition]{Definition}
\aliascntresetthe{definition}

\newtheorem*{definition*}{Definition}
\newtheorem*{notation*}{Notation}
\newtheorem*{organization*}{Organization}
\newtheorem*{theorem*}{Theorem}

\numberwithin{equation}{section}

\usepackage{hyperref}

\hypersetup{
    colorlinks=true,
    linkcolor={blue!90!black},
    citecolor={red!80!black},
    urlcolor={cyan}
}

\usepackage[capitalise,nameinlink,noabbrev]{cleveref}

\crefname{thm}{Theorem}{Theorems}
\crefname{lemma}{Lemma}{Lemmas}
\crefname{cor}{Corollary}{Corollaries}
\crefname{prop}{Proposition}{Propositions}
\crefname{conj}{Conjecture}{Conjectures}
\crefname{question}{Question}{Questions}
\crefname{definition}{Definition}{Definitions}
\crefname{remark}{Remark}{Remarks}
\crefname{example}{Example}{Examples}
\crefname{notation}{Notation}{Notations}

\crefname{innercustomthm}{Theorem}{Theorems}
\crefname{innercustomcor}{Corollary}{Corollaries}
\crefname{innercustomprop}{Proposition}{Propositions}
\crefname{innercustomconj}{Conjecture}{Conjectures}
\crefname{innercustomquestion}{Question}{Questions}

\crefname{equation}{Equation}{Equations}
\crefname{section}{Section}{Sections}
\crefname{subsection}{Subsection}{Subsections}
\crefname{figure}{Figure}{Figures}
\crefname{table}{Table}{Tables}
\crefname{enumi}{Item}{Items}

\newcommand{\Gr}{\mathbf{Gr}}
\newcommand{\kk}{\mathbb{K}}
\newcommand{\bbQ}{\mathbb{Q}}

\newcommand{\bfa}{\mathbf{a}}
\newcommand{\lex}{\mathrm{lex}}

\newcommand{\rad}{\mathrm{rad}}

\newcommand{\GL}{\mathrm{GL}}

\newcommand{\Quot}{\mathrm{Quot}}
\newcommand{\PP}{\mathbf{P}}
\newcommand{\mfh}{\mathfrak{h}}
\newcommand{\mfm}{\mathfrak{m}}
\newcommand{\bfx}{\boldsymbol{x}}
\newcommand{\bff}{\boldsymbol{f}}
\newcommand{\bfu}{\boldsymbol{u}}
\newcommand{\bfv}{\boldsymbol{v}}
\newcommand{\bfe}{\mathbf{e}}

\newcommand{\EE}{\mathscr{E}}
\newcommand{\CC}{\mathscr{C}}
\newcommand{\FF}{\mathscr{F}}
\newcommand{\OO}{\mathscr{O}}
\newcommand{\MM}{\mathscr{M}}
\newcommand{\RR}{\mathscr{R}}
\newcommand{\II}{\mathscr{I}}
\newcommand{\NN}{\mathscr{N}}
\newcommand{\KK}{\mathscr{K}}
\newcommand{\GG}{\mathscr{G}}
\newcommand{\LLL}{\mathscr{L}}
\newcommand{\HH}{\mathcal{H}}
\newcommand{\LL}{\mathcal{L}}

\newcommand{\QQ}{\mathcal{Q}}
\newcommand{\UU}{\mathcal{U}}
\newcommand{\VV}{\mathcal{V}}
\newcommand{\XX}{\mathcal{X}}
\newcommand{\YY}{\mathcal{Y}}
\newcommand{\ZZ}{\mathbf{Z}}

\DeclareMathOperator{\Hom}{Hom}

\DeclareMathOperator{\Ext}{Ext}
\DeclareMathOperator{\End}{End}
\DeclareMathOperator{\Def}{Def}
\DeclareMathOperator{\Fr}{Fr}
\DeclareMathOperator{\Supp}{Supp}
\DeclareMathOperator{\codim}{codim}
\DeclareMathOperator{\Hilb}{Hilb}
\DeclareMathOperator{\Proj}{Proj}
\DeclareMathOperator{\Spec}{Spec}
\DeclareMathOperator{\Flag}{Flag}
\DeclareMathOperator{\Ann}{Ann}
\DeclareMathOperator{\ResQuot}{ResQuot}
\DeclareMathOperator{\Sym}{Sym}

\DeclareMathOperator{\Bl}{Bl}

\begin{document}
\author[R.\,Ramkumar]{Ritvik~Ramkumar}
\address{Department of Mathematics\\University of Notre Dame\\South Bend, IN\\USA}
\email{rramkuma@nd.edu}
\title{Quot schemes on projective space}
\maketitle
\begin{abstract}
We develop general methods for studying Quot schemes on projective space, including a new construction of components with no common neighbours.
We show that the radius of Hilbert and Quot schemes can grow linearly with the dimension of the subschemes being parametrized, and construct the first examples of Hilbert schemes of curves with radius two. 
For absolutely split bundles on projective space, we construct a distinguished generically smooth component and classify all smooth Quot schemes. 
In most cases, we give explicit geometric descriptions in terms of products of generalized flag varieties and projective bundles over Grassmannians.
\end{abstract}

\section{Introduction}
Grothendieck's Quot scheme \cite{Grothendieck__construction_Hilb} is a ubiquitous moduli space in algebraic geometry. 
For example, it underlies constructions of moduli spaces of curves, semistable sheaves, and stable varieties \cite{Gieseker_ModuliCurves,HL,Kollar_Families}, while its
intersection theory and virtual invariants play a central role in enumerative geometry \cite{MarianOprea2007,MNOP}.
Hilbert and Quot schemes also have deep connections to combinatorics and representation theory \cite{Haiman_macdonald,Nakajima1997}.
Their geometry, however, can be very complicated, with multiple irreducible components and nonreduced structure even when parameterizing smooth varieties \cite{Mumford_nonreduced_component}. 
In fact, it is known that that all finite type singularities can occur on Quot schemes \cite{Vakil_MurphyLaw,Jelisiejew__Pathologies}.
Although understanding the global geometry of these schemes, particularly the incidence of their irreducible components, is a natural problem, such pathological behavior has discouraged its systematic study.

Despite this complexity, Hilbert schemes of projective space exhibit remarkable structural properties. 
A key reason is the action of the general linear group, particularly degenerations to Borel-fixed points, which allow one to connect complicated deformation theory to more tractable combinatorics.
A foundational example is Hartshorne's connectedness theorem for $\Hilb^P(\PP^n)$ \cite{har66connectedness}. 
Building on this, Reeves and Stillman identified a distinguished generically smooth component \cite{RS97}, while Reeves showed that every irreducible component lies at a distance of at most $\deg(P)+1$ from it \cite{Reeves_radius}. 
More recently, Skjelnes and Smith classified the smooth Hilbert schemes of projective space, revealing a surprisingly rich collection of examples \cite{Staal_smooth_hilbert_schemes,SS23}.

The corresponding questions for Quot schemes remain considerably less understood. 
Our goal is to develop a systematic framework for studying Quot schemes on projective space, with particular emphasis on absolutely split bundles. 
We study their smoothness as well as the global structure and incidence of irreducible components. 
While our classification results exploit the structure of split bundles, our methods have broader applicability; for instance, we prove that the Quot schemes of the tangent bundle on $\PP^2$ are connected. 
We expect further results in this direction and conclude with several open questions.

We now give an overview of our results.
Let $\kk$ be an algebraically closed field and $S := \kk[x_0,\dots,x_n]$ denote the coordinate ring of $\PP^n$. 
A \textbf{partition} is a weakly decreasing sequence $\lambda=(\lambda_1,\ldots,\lambda_{\ell})$ of positive integers, possibly empty.
For a weakly increasing sequence of integers $\bfa=(a_1,\ldots,a_r)$, consider the absolutely split bundle
$$
\EE_{\bfa}=\bigoplus_{i=1}^r\OO_{\PP^n}(-a_i).
$$
Fix $0\leq s<r$ and $q=r-s$, and let $\lambda$ be a partition with parts of size at most $n$. 
Consider the polynomials
$$
p_\lambda(t)=\sum_{i=1}^{\ell} \binom{t+\lambda_i-i}{\lambda_i-1},
\qquad
P_{\bfa,s,\lambda}(t) = \sum_{j=q+1}^r\binom{t-a_j+n}{n}+p_\lambda(t-a_q),
$$
and let
$$
\QQ_{\bfa}(s,\lambda) = \Quot^{P_{\bfa,s,\lambda}}_{\PP^n}(\EE_{\bfa}).
$$
The number $s$ is the rank of the quotient being parameterized and $q$ is the rank of its kernel.
By explicitly exhibiting a submodule, called the \textbf{lexicographic module} (see \Cref{lexmodule}), it follows from Macaulay's theorem that the schemes $\QQ_{\bfa}(s,\lambda)$ are precisely the nonempty Quot schemes associated to an absolutely split source bundle \cite{Hulett95, Pardue} .

\subsection{General constructions on projective space and Borel-fixed points}
We prove two general results that underpin much of our approach throughout the paper.
The first compares deformations of a quotient $[\EE\twoheadrightarrow\FF]$, with kernel $\KK$, to deformations of the quotient $[\KK^{\vee\vee}\twoheadrightarrow\KK^{\vee\vee}/\KK]$; under the hypotheses of \Cref{double_dual_prop}, the corresponding points on the two Quot schemes are smoothly equivalent. 
The second concerns the behavior of Quot schemes under replacement of the source by a direct summand. 
Both results are formulated for vector bundles on an arbitrary smooth projective scheme, and provide a flexible framework that we expect to be useful beyond the absolutely split case considered here. 
Their underlying ideas recur throughout the paper, playing a central role in our study of smoothness and the geometry of irreducible components.

Pardue proved that $\QQ_{\bfa}(s,\lambda)$ is connected \cite{Pardue}, and our first result establishes the existence of a distinguished generically smooth component and a criterion for when an arbitrary Borel-fixed point lies on this component.
Note that connectedness fails for arbitrary source bundles, as we exhibit a rank-two bundle $\EE$ on $\PP^2$ for which $\Quot_{\PP^2}^{t+1}(\EE)$ consists of two reduced points (see \Cref{disconnected_ex}). Thus, this is essentially the best possible statement in such generality.

\begin{customthm}{A}
For $n\geq 2$, the lexicographic point of $\QQ_{\bfa}(s,\lambda)$ is smooth and hence belongs to a unique irreducible component, called the
\textbf{lexicographic component}. Moreover, a Borel-fixed point lies on the lexicographic component if and only if it has the same double saturation as the lexicographic point.
\end{customthm}

See \Cref{double_sat} for the definition of double saturation. 
This result is analogous to the results of Reeves--Stillman and Reeves for the Hilbert scheme \cite{RS97, Reeves_radius}.

Rather surprisingly, the lexicographic point can be singular when $n=1$, and \Cref{non_distinguish} shows there is no way to pick a distinguished component just based on the Borel-fixed points. 
Nevertheless, there is a component that contains all the Borel-fixed points.

\begin{customthm}{B}
For $n=1$, all the Borel-fixed points of $\QQ_{\bfa}(s,\lambda)$ lie on a common irreducible component.
\end{customthm}

\subsection{Smooth Quot schemes}
We now describe the lexicographic component and classify the smooth Quot schemes. 
The geometric constructions are built on the residual flags of Skjelnes--Smith \cite{SS23}, and give a modular interpretation of the lexicographic component. 
There are two stages; for $n\geq 2$, Quot schemes corresponding to partitions with no parts equal to one admit descriptions in terms of generalized flag varieties and projective bundles. 
We refer to these as the \textbf{base Quot schemes}. 
The remaining parts equal to one are then incorporated through successive quotients of their universal kernels.

Consider a non-empty partition $\mu = (n^c,\nu)$, where $\nu$ has no parts equal to $n$ or $1$. 
Let $q_0$ denote the number of occurences of $a_q$ among the first $q$ entries of $\bfa$ and $s_0$ the number of occurences of $a_q$ among the last $s$ entries of $\bfa$.
Let $\mathbf{G}=\Gr(s_0+1,q_0+s_0)$ and  $\mathscr{Q}$ denote the universal quotient bundle on $\mathbf{G}$. 
We define the locally free sheaf $\mathscr V_c$ on $\mathbf{G}$ by
$$
\mathscr V_c=
(S_c\otimes_\kk \mathscr{Q}) \oplus \bigoplus_{a_j>a_q} S_{a_q+c-a_j}\otimes_\kk\OO_{\mathbf G}.
$$

\begin{customthm}{C} 
Let $n\geq 2$. 
For a nonempty partition $\mu = (n^c,\nu)$, where $\nu$ has no parts equal to $n$ or $1$, we have
$$
\QQ_{\bfa}(s,\mu) \simeq \Flag(\nu;\PP^n)\times\PP_{\mathbf G}(\mathscr V_c).
$$
\end{customthm}
We now explain how to pass to the general Quot scheme. 
Let $\lambda = (\mu,1^b)$ with $\mu$ having no parts equal to one.
Let $\mathscr{K}_\mu$ denote the universal kernel on $\PP^n\times\QQ_{\bfa}(s,\mu)$.
Define the \textbf{residual Quot scheme}
$$
\ResQuot_{\bfa}(s,\lambda)
:= \Quot^b_{\PP^n\times\QQ_{\bfa}(s,\mu)/ \QQ_{\bfa}(s,\mu)}(\mathscr K_\mu).
$$
It parameterizes residual quotient flags $ [\EE_{\bfa}\twoheadrightarrow\FF  \twoheadrightarrow\FF_0]$,
where $[\EE_{\bfa} \twoheadrightarrow \FF_0]$ belongs to $\QQ_{\bfa}(s,\mu)$ and $\ker(\FF \to \FF_0)$ has length $b$.
The locus inside $\ResQuot_{\bfa}(s,\lambda)$ where the kernel is supported at $b$ distinct points of the locally free locus of $\FF_0$ is called the \textbf{configuration locus} (see \Cref{defn_residual} for more details).

\begin{customthm}{D}
Let $n\geq2$ and write $\lambda=(\mu,1^b)$, where $\mu$ has no parts equal to one. 
Forgetting the second quotient defines a natural closed immersion
$$
\ResQuot_{\bfa}(s,\lambda) \hookrightarrow \QQ_{\bfa}(s,\lambda).
$$
The closure of the configuration locus is precisely the lexicographic component.
\end{customthm}

Since the singular locus is closed and Borel-invariant, Borel's fixed point theorem implies that $\QQ_{\bfa}(s,\lambda)$ is smooth if and only if it is smooth at every Borel-fixed point.
Using our description of the Borel-fixed points on the lexicographic component, we are able to classify all the smooth Quot schemes.

\begin{customthm}{E}
The Quot scheme $\QQ_{\bfa}(s,\lambda)$ is smooth if and only if one of the following holds.
\begin{enumerate}
\item $n=1$, $\lambda=(1^b)$, and either $b=0$, $q=1$, or there is no $j>q$ such that $a_1+2\leq a_j\leq a_q+b$.
\item $n=2$ and $q=1$.
\item $n\geq3$, $q=1$, and, writing $\lambda=(n^c,\nu)$ with $\nu$ having no parts equal to $n$, and either
\begin{enumerate}
\item[(a)] $\nu$ has no parts equal to one, or
\item[(b)] $\nu=(1^d)$ with $1\leq d\leq3$, or
\item[(c)] $\nu=(m^d,1)$ with $2\leq m<n$ and $d\geq1$, excluding $(m,d)=(2,2)$ and $(2,3)$.
\end{enumerate}
\item $n\geq2$, $q\geq2$, and either $\lambda= (1)$ or $\lambda$ has no parts equal to one.
\item $n\geq2$, $q\geq2$, and $\lambda=(n^c,1)$ with $c\geq1$,
where either $n\geq3$, or $n=2$ and $c+a_q-a_{q-1}\geq3$.
\end{enumerate}
\end{customthm}
In each case, the Quot scheme is irreducible and, when $n\geq 2$, it is isomorphic to the residual Quot scheme. 
Moreover, for $n\geq2$, the smooth Quot schemes admit explicit geometric descriptions, in most cases as products of generalized flag varieties and projective bundles, with a few exceptional cases described by blow-ups (see \Cref{smooth_quot_geometry}). 
This fully generalizes the recent work of Skeljnes-Smith-Stillman \cite{SSS26}, which treats the case where $\bfa$ is a constant sequence.

\subsection{Radius of Quot schemes}
We now turn to the configuration of irreducible components of Quot schemes, about which very little is known even in the case of Hilbert schemes.
The \textbf{incidence graph} of a connected Noetherian scheme $\XX$ is the graph whose vertices are the irreducible components of $\XX$, with two distinct vertices joined by an edge if the
corresponding components intersect. 
We define the \textbf{radius} of $\XX$, denoted by $\mathrm{rad}(\XX)$, to be the radius of its incidence graph. This is the minimum, over all vertices, of the maximum distance to any other vertex (see \Cref{radius_defn} for more details).

Reeves proved that the radius of $\Hilb^P(\PP^n)$ is at most $\deg(P)+1$ \cite{Reeves_radius}, i.e., it is at most one more than the dimension of the subschemes being parameterized. 
Pardue obtained a closely related bound for $\Quot_{\PP^n}^{P}(\EE_\bfa)$, namely that every closed point can be joined to the lexicographic point by a chain of rational curves of
length at most $\deg(P)+2$ \cite{Pardue}. 
However, without a distinguished component, this does not yield the analogous bound on the radius. 
As we have already seen, the case of $\PP^1$ presents some unexpected phenomena. 
Nonetheless, we are able to show that

\begin{customthm}{F} 
The radius of $\QQ_{\bfa}(s,\lambda)$ is at most $\deg(P_{\bfa,s,\lambda})+1$.
\end{customthm}

Lower bounds for the radius are a more subtle question. 
For constant Hilbert polynomials, the radius is either zero or one, and distinguishing
these cases amounts to determining the irreducibility of a Quot scheme of points, a problem studied extensively (see \cite{iarrobino_reducibility, cartwright_erman_velasco_viray_Hilb8,DJNT,jelisiejew_sivic} and the references therein).
For non-constant polynomials, there has been no progress in establishing lower bounds.\footnote{In fact, an unpublished note of ours from 2020 was the first example of any Hilbert scheme having radius greater than one (\url{https://arxiv.org/abs/2008.11968}).}

We show that the Hilbert scheme of curves can have radius two.
We introduce the notion of a \textbf{TNT pair}, which is a pair of Artinian algebras with trivial negative tangents whose cones have the same Hilbert polynomial (see \Cref{TNT_pair}).
When the Artinian algebras are smooth on their Hilbert scheme, we show that a TNT pair gives rise to generically smooth components parameterizing cones.
A cohomological criterion then shows that these components are disjoint and have no common neighbour (\Cref{TNT_pair_cor}). 
Using this construction, we give the first example of Hilbert schemes of curves of radius two, proving that Reeves' bound is sharp in dimension one; another example of lower degree and genus is given in \Cref{radius_two_14}.

\begin{customthm}{G}  
The Hilbert scheme $\Hilb^{18t-9}(\PP^9)$ has radius two.
It contains the following generically smooth components:
\begin{enumerate}
\item $\HH_1$, whose general member is a cone over an Artinian
algebra with Hilbert function $(1,7,10)$ in a varying $\PP^8$;
\item $\HH_2$, whose general member is a cone over a Gorenstein
algebra having Hilbert function $(1,8,8,1)$;
\item $\mathcal A$, whose general member is a smooth genus-$10$
curve in a varying $\PP^8$;
\item $\mathcal B$, whose general member is a smooth canonical
genus-$10$ curve in $\PP^9$.
\end{enumerate}
A shortest path from $\HH_1$ to $\HH_2$ is as follows
\begin{center}
\begin{tikzpicture}[every node/.style={inner sep=4pt}]
\node (h1) at (0,0) {$\HH_1$};
\node (a) at (2,0) {$\mathcal A$};
\node (b) at (4,0) {$\mathcal B$};
\node (h2) at (6,0) {$\HH_2$};
\draw (h1)--(a)--(b)--(h2);
\end{tikzpicture}
\end{center}

\end{customthm}

A folklore expectation was that the radius of Hilbert schemes would remain bounded by a small constant, far below the general upper bound. 
We show instead that the radius can grow linearly with the dimension of the subschemes being parametrized. 
We prove a general slicing criterion (\Cref{slicing_criterion}) which roughly states that if the hyperplane section of an irreducible component has positive distance from the lexicographic component, then the original component has strictly greater distance.
Applying this repeatedly to projective cones shows that the linear order of growth in Reeves' bound is optimal.

\begin{customthm}{H}
Let $P_d(t)= \binom{t+d}{d} +4\binom{t+d-1}{d} +3\binom{t+d-2}{d}$.
Then the Hilbert scheme $\Hilb^{P_d}(\PP^{d+4})$ contains cones over algebras of length $8$ and satisfies $$
\left\lceil\frac{d+1}{2}\right\rceil \leq \rad\bigl(\Hilb^{P_d}(\PP^{d+4})\bigr) \leq d+1.
$$
\end{customthm}
We decided to restrict our explicit examples to Hilbert schemes because the results are already new
in this setting. 
However, the cone constructions in \Cref{radius_sec} apply more generally to suitable pairs of finite modules with TNT, producing components of Quot schemes at large distance from one another in arbitrarily large projective spaces. 
The study of TNT modules is an active area of research \cite{jelisiejew_sivic, Jelisiejew__Elementary,Satriano_Staal, SS24} and offers a source of further examples. For instance, allowing source bundles of higher rank gives examples of smaller degree than those in \Cref{radius_two_14}.

\subsection{Conjectures and questions}
The geometry of Hilbert and Quot schemes on other rational varieties, such as Grassmannians and toric varieties, is much less understood. 
Recently, Cid-Ruiz \cite{CidRuiz26} showed that the multigraded Hilbert schemes of $\PP^2 \times \PP^1$ can be disconnected. 
On the other hand, we prove in \Cref{connected_tan} that every nonempty Quot scheme of the tangent bundle
of $\PP^2$ is connected. 
This suggests that, rather than changing the base space, fixing it to be $\PP^n$ and varying the source bundle may lead to more general connectedness results.

\begin{question}
Determine when $\Quot_{\PP^n}^P(\EE)$ is connected, for a locally free sheaf $\EE$ on $\PP^n$. 
\end{question}
We expect that there are many source bundles for which every nonempty Quot scheme is connected. One can then ask further questions about non-emptiness, the number of irreducible components, and so on. 
A natural starting point would be stable bundles on $\PP^n$.

Although Quot schemes on $\PP^1$ can be singular and reducible, Popa and Roth \cite[Theorem~6.4]{PopaRoth} proved that they are irreducible and generically smooth in sufficiently large degree.
However, we expect that they are always reduced, regardless of the degree. 

\begin{conj}
Let $\EE$ be a locally free sheaf on $\PP^1$.  Then $\Quot_{\PP^1}(\EE)$ is reduced.
\end{conj}

Another natural direction is to study the motivic invariants of these Quot schemes. 
Explicit formulas are known for the motivic classes of Quot schemes of points on smooth curves and surfaces \cite{BFP20,Moz19}, and for positive-rank quotients of split bundles on curves under smoothness and unobstructedness assumptions \cite{MR25}.
Our geometric descriptions of smooth Quot schemes provide a starting point for similar computations on projective space.
It would be interesting to understand how much of this structure persists beyond the smooth cases.

\begin{question}
Compute the class $[\QQ_{\bfa}(s,\lambda)]$ in the Grothendieck ring of varieties.
\end{question}

Finally, we still have the question of whether Reeves' bound is sharp in all dimensions.
 
\begin{question} Determine whether the bound $\rad(\QQ_{\bfa}(s,\lambda)) \leq \deg(P_{\bfa,s,\lambda})+1$ is sharp in all dimensions.
\end{question}
We believe that the methods of \Cref{radius_sec} can be extended to construct a Hilbert scheme of surfaces with radius three. 
For instance, with additional cohomological restrictions, it may be possible to arrange suitable cone constructions in a controlled combinatorial pattern to obtain the desired bound.

\section*{Acknowledgments}
I would like to thank Paolo Lella for reigniting my interest in the radius. 
I would also like to thank Joachim Jelisiejew and Greg Smith  for helpful comments.
Many of the ideas in this paper grew out of questions I have considered over the years but never committed to writing. With artificial intelligence playing an increasingly prominent role in mathematical research, I felt that this was a good time to put these ideas on record, in the hope that they might stimulate further discussion and lead to new developments.

\section{Deformation theory of Quot schemes}
In this section, we establish some general deformation-theoretic results for Quot schemes and study their behavior under restriction to locally free subsheaves. 

\begin{notation} Unless otherwise specified, all schemes are of finite type over an algebraically closed field $\kk$. Given two $\kk$-schemes $X$ and $T$ and a coherent sheaf $\FF$ on $X$, we write $X_T = X \times _\kk T$ and $\FF_T$ for the pullback of $\FF$ to $X_T$. 
We use $\PP(V)$ for the projective space of lines in $V$ and $\Gr(k,V)$ for the Grassmannian of rank $k$ quotients of $V$, with the analogous convention for vector bundles.
A projective scheme always means a projective $\kk$-scheme. 
If $X$ is projective and $\OO_X(1)$ is an ample line bundle, the Hilbert polynomial of $\FF$ is $P_{\FF}(t) = \chi(\FF(t))$.

By a \textbf{small extension} with kernel $I$ we mean a surjection $A' \twoheadrightarrow A$ of local Artinian $\kk$-algebras with kernel $I$ satisfying $\mfm_{A'}I=0$.
\end{notation}

\begin{definition}[{\cite[Theorem 2.2.4]{HL}}]
Let $X$ be a projective $\kk$-scheme with a fixed ample invertible sheaf $\OO_X(1)$, let $\EE$ be a coherent sheaf on $X$, and let $P\in \bbQ[t]$ be a numerical polynomial. 
Consider the contravariant functor on finite type $\kk$-schemes defined by
$$
T \mapsto
\left\{\EE_T\twoheadrightarrow\FF \ \middle|\
\begin{array}{l}
\FF\text{ is coherent and flat over } T,\\
P_{\FF_t}=P\text{ for every } t \in  T
\end{array}
\right\}
\big/\sim,
$$
where Hilbert polynomials are computed with respect to $\OO_X(1)$, and two quotients are equivalent if their kernels coincide. 
The \textbf{Quot scheme}, denoted by $\Quot_X^P(\EE)$, is the projective $\kk$-scheme representing this functor. 
We write $\Quot_X(\EE)=\bigsqcup_P\Quot_X^P(\EE)$ for the full Quot scheme.
\end{definition}

\begin{notation}
Closed points of $\Quot_X(\EE)$ will be denoted by $[\EE \twoheadrightarrow \FF]$. 
As long as there is no ambiguity, we will refer to any such point as a quotient.
Its kernel is $\ker(\EE \twoheadrightarrow \FF)$.
\end{notation}

\begin{definition}[{\cite[Section 2.A.1]{HL}}]
Let $P_1,P_2 \in \bbQ[t]$ be numerical polynomials. The \textbf{flag Quot scheme}, denoted by $\mathrm{FlagQuot}_X^{P_1,P_2}(\EE)$, is the projective $\kk$-scheme representing the functor on finite type $\kk$-schemes defined by
$$
T \mapsto
\left\{\EE_T \twoheadrightarrow\FF_1\twoheadrightarrow\FF_2 \ \middle|\
\begin{array}{l}
\FF_i\text{ is coherent and flat over } T,\\
P_{(\FF_i)_t}=P_i\text{ for every } t \in T \text{ and } i
\end{array}
\right\}
\big/\sim,
$$
where two flags are equivalent if the kernels of the maps
$\EE_T\to\FF_i$ coincide for both $i$. 
Forgetting the second quotient induces a morphism
$$
\mathrm{FlagQuot}_X^{P_1,P_2}(\EE) \longrightarrow\Quot_X^{P_1}(\EE),
\qquad
[\EE\twoheadrightarrow\FF_1\twoheadrightarrow\FF_2] \mapsto [\EE\twoheadrightarrow\FF_1].
$$
\end{definition}

\begin{lemma}[{\cite[Proposition 2.2.7]{HL}}]
The tangent space to the closed point $[\EE\twoheadrightarrow\FF]\in\Quot_X^P(\EE)$ with
kernel $\KK$ is $T_{[\EE\twoheadrightarrow\FF]}\Quot_X^P(\EE)\cong\Hom_X(\KK,\FF)$.
\end{lemma}

Our first goal is to understand the locus of quotients that remain surjective when restricted to a locally free subsheaf.

\begin{definition} 
Let $\EE$ be a locally free sheaf on $X$ and $\EE' \subseteq \EE$ a locally free subsheaf. 
The locus in $ \Quot_X(\EE)$ where the restriction to $\EE'$ is surjective is defined to be
$$
\UU_{\EE'} = \{[\EE \twoheadrightarrow \FF] \in \Quot_X(\EE): \EE' \hookrightarrow \EE \to \FF \text{ is surjective}\}
$$ 
\end{definition}

\begin{lemma} \label{restriction_open}
Let $X$ be a projective scheme and $\EE$ a locally free sheaf. Let $\EE'$ be a locally free subsheaf with cokernel $\mathscr{C} = \EE/\EE'$. 
Then  $\UU_{\EE'}$ is open and restriction induces an affine morphism
$$
\Psi:\UU_{\EE'} \longrightarrow \Quot(\EE'),
$$
whose nonempty fiber over $[\EE' \twoheadrightarrow \FF]$ is a torsor under the additive group $\Hom_X(\mathscr{C},\FF)$.
\end{lemma}
\begin{proof}
Let $\QQ=\Quot_X(\EE)$ and $\QQ' = \Quot_X(\EE')$, and let $q: \EE_\QQ \twoheadrightarrow \MM$ and $q': \EE'_{\QQ'} \twoheadrightarrow \MM'$ be the respective universal quotients. 
Let $\pi: X \times \QQ \to \QQ$ be the projection and let $\GG = \mathrm{coker}(q|_{\EE'_\QQ})$. 
The locus $\UU_{\EE'} = \QQ \setminus \pi(\Supp(\GG))$ is open since $\Supp(\GG)$ is closed and $\pi$ is proper. 
Thus, we get a natural map $\Psi:\UU_{\EE'} \to \Quot(\EE')$ which, explicitly at the level of $T$-points, corresponds to restricting the surjection $[\EE_T \twoheadrightarrow\FF]$ to $[\EE'_T \twoheadrightarrow \FF]$.

By \cite[\href{https://stacks.math.columbia.edu/tag/08K6}{Tag 08K6}]{stacks_project}, the relative Hom functors $\underline{\mathrm{Hom}}_{X\times \QQ'/\QQ'}(\EE_{\QQ'},\MM')$ and $\underline{\mathrm{Hom}}_{X\times \QQ'/\QQ'}(\EE'_{\QQ'}, \MM')$ are represented by affine $\QQ'$-schemes $\HH$ and $\HH'$, respectively. 
Now we have a natural restriction morphism $\HH \to \HH'$ and a morphism $\QQ' \to \HH'$ given by the universal quotient. 
It is straightforward to check that $\UU_{\EE'} \simeq \HH \times_{\HH'} \QQ'$ and that the morphism $\Psi$ is just the projection onto the second component. 
Since $\HH$ and $\HH'$ are affine over $\QQ'$, the restriction $\HH\to\HH'$, and thus $\Psi$, are affine.

Finally, every nonempty fiber is a fiber of the linear restriction map $\Hom_X(\EE, \FF) \to \Hom_X(\EE',\FF)$, and hence a torsor under its kernel $\Hom_X(\CC,\FF)$.
\end{proof}

\begin{definition} 
An element of the fiber $\Psi^{-1}([\EE' \twoheadrightarrow \FF])$ is said to \textbf{extend} $[\EE' \twoheadrightarrow \FF]$.
\end{definition}

\begin{lemma} \label{fibration}
Let $X$ be a projective scheme and $\EE$ a locally free sheaf.
Let
$$
0\to\mathscr E'\to\mathscr E\to\mathscr C\to0
$$
be an exact sequence of locally free sheaves on $X$ and assume that $\Quot_X(\EE')$ is smooth at $z=[\EE' \twoheadrightarrow \FF]$.
\begin{enumerate}
\item If $\Ext^1_X(\CC,\FF) = 0$, then $\Quot_X(\EE)$ is smooth at every point extending $z$.
\item If $\Psi$ admits a section $\sigma$ over a neighbourhood of $z$, then $\Quot_X(\EE)$ is smooth at $\sigma(z)$ if and only if the function $h^0(X, \CC^\vee \otimes -)$ is locally constant in a neighbourhood of $z$.
\end{enumerate}
\end{lemma}
\begin{proof}
Let $\QQ'=\Quot_X^P(\EE')$ and let $\EE'_{\QQ'}\twoheadrightarrow  \MM'$ be the universal quotient, with kernel $\KK'$.  
Consider the exact sequence
$$
0\longrightarrow  \EE'_{\QQ'}/\KK' \longrightarrow \EE_{\QQ'}/\KK' \longrightarrow  \CC_{\QQ'} \longrightarrow 0.
$$
By \cite[Theorem 2.2.4]{HL}, we can consider the relative Quot scheme $\Quot_{X\times \QQ'/\QQ'}^P(\EE_{\QQ'}/\KK')$, which is projective over $\QQ'$. 
Note that the morphism $\Psi$ identifies $\mathcal U_{\EE'}$ with the following locus in the relative Quot scheme
$$
\{[\EE_{\QQ'}/\KK' \twoheadrightarrow \mathscr{H}]: \EE'_{\QQ'}/\KK' \to \EE_{\QQ'}/\KK' \to  \mathscr{H} \text{ is an isomorphism} \}.
$$
The same argument as in the previous lemma shows that this locus is open.

Assume that we are in case (i) and consider any point $y$ extending $z$. 
Since $\QQ'$ is smooth at $z$, it suffices to show that $y$, which we may assume is of the form $[\EE_{\QQ'}/\KK' \twoheadrightarrow \FF]$, is smooth in the relative Quot scheme.
The kernel of the map $\EE_{\QQ'}/\KK' \twoheadrightarrow \FF$ is exactly $\CC$.
Thus, by \cite[Lemma 2.2.6]{HL}, the obstruction to lifting $[\EE_{\QQ'}/\KK' \twoheadrightarrow \FF]$ over a small extension with kernel $I$ lies in $\Ext_X^1(\mathscr C,\FF)\otimes_{\kk} I$.
By assumption, this group vanishes and the infinitesimal criterion for smoothness implies that the relative Quot morphism is smooth at $y=[\EE_{\QQ'}/\KK'\twoheadrightarrow\FF]$, as required.

Now assume that we are in case (ii). 
By shrinking, we may assume that $\sigma$ is defined over a smooth connected open neighbourhood $\VV$ of $z$. 
Let $d = \dim(\VV)$, and for any point $v=[\EE'\twoheadrightarrow\MM_v] \in \VV$, write $h(v)=\dim_{\kappa(v)}\Hom_X(\CC_v,\MM_v)$.
Since $\Psi\circ\sigma=\operatorname{id}_\VV$, the differential of $\Psi$ is surjective at $\sigma(v)$. 
By \Cref{restriction_open}, its kernel is the tangent space to an affine space of dimension $h(v)$.
In particular, for every closed point $v\in \VV$, we have  $\dim T_{\sigma(v)}\Quot_X(\EE)=d+h(v)$.

Suppose first that $h$ is constant and, by abuse of notation, let its constant value also be $h$.
By cohomology and base change, the relative Hom functor is a vector bundle of rank $h$.
The morphism $\Psi^{-1}(\VV)\to\VV$ is a torsor under this bundle. 
The section $\sigma$ trivializes this torsor, so $\Psi^{-1}(\VV)$ is the total space of a rank-$h$ locally free sheaf over $\VV$. In particular, it is smooth of dimension $d+h$. 
Since $\dim T_{\sigma(v)}\Quot_X(\EE)=d+h$ for every closed point $v\in\VV$, the Quot scheme is smooth at $\sigma(v)$.

Conversely, suppose that $\Quot_X(\EE)$ is smooth at $\sigma(z)$.
After shrinking $\VV$, its image under $\sigma$ lies in a smooth open subset of $\Quot_X(\EE)$ of constant dimension $D$.
The tangent-space formula above then gives $h(v)=D-d$ for every closed point $v\in \VV$, as required. 
\end{proof}

We now prove a proposition which can be thought of as a vast generalization of the fact that one can detach hypersurfaces. 

\begin{definition}[{\cite{Jelisiejew__Pathologies, Vakil_MurphyLaw}}] Two pointed schemes $(\XX,x)$ and $(\YY,y)$ are \textbf{smoothly equivalent} if there exists a third pointed scheme $(\mathcal{Z},z)$ with smooth maps $(\mathcal{Z},z) \to (\XX,x)$ and $(\mathcal{Z},z) \to (\YY,y)$.
\end{definition}

\begin{thm} \label{double_dual_prop}
Let $X$ be a smooth projective scheme of dimension at least two, and $\EE$ a locally free sheaf. Let $z = [\EE \twoheadrightarrow \FF]$ be a point on $\Quot_X(\EE)$ with kernel $\KK$ and assume that either
\begin{enumerate}
\item[(i)] $\KK$ has rank one and $H^1(X,\OO_X)=0$, or
\item[(ii)] $\Ext^1_X(\KK,\KK^{\vee\vee}) = 0$.
\end{enumerate}
Then the points $[\EE \twoheadrightarrow \FF]$ and $[\KK^{\vee\vee} \twoheadrightarrow \KK^{\vee\vee}/\KK]$, on their respective Quot schemes, are smoothly equivalent.
 \end{thm}
\begin{proof} 
For simplicity, let $z = [\EE \twoheadrightarrow \FF]$, $\RR=\KK^{\vee\vee}$ and $z'=[\mathscr R\twoheadrightarrow\mathscr R/\KK]$.
Taking double duals of the inclusion $i:\KK\hookrightarrow\EE$ gives $f_0=i^{\vee\vee}:\RR\to\EE^{\vee\vee}=\EE$.
By the universal property of the reflexive hull \cite[\href{https://stacks.math.columbia.edu/tag/0EBH}{Tag 0EBH}]{stacks_project}, this is the unique extension of $i$.
It is injective because it is generically injective and $\mathscr R$ is torsion-free.
Let $j:U\hookrightarrow X$ be the inclusion of the open locus where $\KK$ is locally free.
Since $\KK$ is torsion-free and $X$ is normal, $X\setminus U$ has codimension at least two, and the canonical map $\KK|_U\to\mathscr R|_U$ is an isomorphism.

For a local Artinian $\kk$-algebra $A$ with residue field $\kk$, let $X_A=X\times\Spec A$ and let $j_A: U_A\hookrightarrow X_A$ be the inclusion.
Since $\RR$ and $\EE$ are reflexive, Hartogs theorem \cite[\href{https://stacks.math.columbia.edu/tag/0EBJ}{Tag 0EBJ}]{stacks_project} gives $\RR \cong j_\star(\RR|_U)$ and $\EE \cong j_\star(\EE|_U)$.
Using flat base change we obtain $\RR_A \cong (j_{A})_\star(\RR_A|_{U_A})$ and  $\EE_A\cong (j_{A})_\star(\EE_A|_{U_A})$ .

Let $\operatorname{Def}_z$ and $\operatorname{Def}_{z'}$ denote the local deformation functors of the respective Quot schemes. 
Consider the functor $\mathrm{D}$ that associates to an Artinian $\kk$-algebra $A$ the set
$$
\mathrm{D}(A)= \left\{ (\MM_A,f_A): [\RR_A/\MM_A] \in \Def_{z'}(A) \text{ and } f_A:\RR_A \to \EE_A \text{ such that } f_A\otimes_A \kk = f_0 \right\}.
$$
There are natural transformations
$$ \Def_{z'}  \xleftarrow{\ \pi } \mathrm{D} \xrightarrow{\ \Xi }\Def_z$$
defined by
$\pi_A(\MM_A,f_A) =[\RR_A \twoheadrightarrow \RR_A/ \MM_A]$ and $\Xi_A(\MM_A,f_A)\mapsto [\EE_A\twoheadrightarrow\EE_A/f_A(\MM_A)]$.
We will show that both $\pi$ and $\Xi$ are well defined and formally smooth. 

Clearly, $\pi$ is well defined and formally smooth. 
Since $\RR_A/\MM_A$ is $A$-flat, the sheaf $\MM_A$ is $A$-flat and its closed fiber is $\KK$. 
The reductions of $f_A$ and $f_A|_{\MM_A}$ modulo the maximal ideal of $A$ are the injections $f_0:\RR\hookrightarrow\EE$ and $\KK\hookrightarrow\EE$, respectively.
The fiberwise injectivity criterion implies that both maps are injective and have $A$-flat cokernels \cite[\href{https://stacks.math.columbia.edu/tag/00ME}{Tag 00ME}]{stacks_project}; thus $\Xi_A$ is well defined. 

Consider the functor $\mathrm{G}$ of infinitesimal automorphisms of $\RR$, namely for every local Artinian $\kk$-algebra $A$ we have
$$
\mathrm{G}(A)= \ker(\mathrm{Aut}_{X_A}(\RR_A) \to \mathrm{Aut}_X(\RR)).
$$
Note that we have
$$
G(A)=1+\End_X(\RR) \otimes_\kk \mfm_A.
$$
Indeed, if $g\in G(A)$, then $g-\mathrm{id}_{\RR_A}$ is zero on the closed fiber, so it belongs to $\End_X(\RR)\otimes_\kk\mfm_A$.
Conversely, for every $u\in\End_X(\RR)\otimes_\kk\mfm_A$, the map $1+u$ is invertible since $u$ is nilpotent; it also restricts to the identity on the closed fiber.
Thus $G(A)=1+\End_X(\RR)\otimes_\kk\mfm_A$.
In particular, this implies that $G(A')\to G(A)$ is surjective for every surjection $A' \twoheadrightarrow A$ of local Artinian $\kk$-algebras.

The group $\mathrm{G}(A)$ acts on $\mathrm{D}(A)$ by $g\cdot(\MM_A,f_A) =(g(\MM_A),f_Ag^{-1})$. Moreover, the map $\Xi_A$ is invariant under this action. 
We will now show that $\Xi$ induces an isomorphism $\mathrm{D}/\mathrm{G} \simeq \Def_z$.  
Let $A$ be a fixed local Artinian $\kk$-algebra.

To show injectivity, assume that $f_A(\MM_A)=f'_A(\MM'_A)$ inside $\EE_A$.
Since $\RR/\KK$ vanishes on $U$, Nakayama's lemma implies that both $\RR_A/\MM_A$ and $\RR_A/\MM'_A$ vanish on $U_A$. 
Hence
$$
\MM_A|_{U_A}=\MM'_A|_{U_A}=\RR_A|_{U_A}.
$$
In particular, the two injective maps $f_A$ and $f'_A$ have the same image on $U_A$.
This allows us to find a unique automorphism $g_{U_A}$ of $\RR_A|_{U_A}$ for which $f'_A g_{U_A}=f_A$.
By Hartogs' theorem and flat base change, $g_{U_A}$ extends uniquely to an automorphism $g$ of $\RR_A$.
We have $f'_A g=f_A$ globally and the restriction to the closed fiber is $f_0\overline{g} = f_0$.
Since $f_0$ is injective, $\overline{g}=\mathrm{id}$, and thus $g\in G(A)$.
Since $f'_A(\MM'_A) = f_A(\MM_A) = f'_Ag(\MM_A)$, the injectivity of $f_A'$ implies that $\MM'_A = g(\MM_A)$. 
This proves the desired injectivity.

For surjectivity, consider an element $[\EE_A \twoheadrightarrow \FF_A] \in \Def_z(A)$. It corresponds to 
an exact sequence
$$
0\longrightarrow \KK_A \xrightarrow{\kappa_A} \EE_A \longrightarrow\FF_A \longrightarrow 0
$$
We will first use the hypothesis in (i) or (ii) to construct a map $i_A: \KK_A \to \RR_A$ which restricts to the  inclusion $\KK\hookrightarrow\RR$.
Once this is done, we will take $\MM_A=i_A(\KK_A)$ to obtain the desired element of $\mathrm{D}(A)$.

Given the hypothesis in (i), $\RR$ is a line bundle.
Since $\KK_A$ is $A$-flat and $X_A\to\Spec A$ is smooth, the sheaf $\KK_A$ is perfect and hence admits a determinant line bundle $\det(\KK_A)$ \cite[Chapter I, Theorem 1]{KM76}.
Its closed fiber is $\det(\KK)\simeq\RR$, and its restriction to $U_A$ is $\KK_A|_{U_A}$.
The Picard exact sequence \cite[\href{https://stacks.math.columbia.edu/tag/0C6R}{Tag 0C6R}]{stacks_project} associated to a small extension $A\twoheadrightarrow A_0$ with kernel $I$ gives
$$
H^1(X,\OO_X)\otimes_\kk I
\longrightarrow \operatorname{Pic}(X_A)
\longrightarrow \operatorname{Pic}(X_{A_0}).
$$
Since $H^1(X,\OO_X)=0$, the restriction map on Picard groups is injective.
Applying this successively along a sequence of small extensions from $A\to \kk$, we see that
$\operatorname{Pic}(X_A)\to\operatorname{Pic}(X)$ is injective.
Since $\det(\KK_A)$ and $\RR_A$ both restrict to $\RR$ on $X$, they are isomorphic.
On $U_A$, the sheaf $\KK_A$ is a line bundle and hence isomorphic to its determinant. 
Putting everything together we obtain
$$
i_A: \KK_A\longrightarrow (j_{A})_\star(\KK_A|_{U_A}) \cong (j_{A})_\star(\det(\KK_A)|_{U_A})
\cong\det(\KK_A)
\cong \RR_A.
$$
Finally, $i_A \otimes \kk$ and the inclusion $\KK \hookrightarrow \RR$ agree on the open set $U$. 
Since $\RR$ is torsion-free, they agree everywhere.

In case (ii), we construct $i_A$ by induction on the length of $A$. Assume $A \twoheadrightarrow B$
is a small extension with kernel $J$ and that $i_B$ has been constructed.
Tensoring $0 \to J \to A \to B \to 0$ with $\RR_A$ and using the fact that $\mfm_AJ = 0$ we obtain
$ 0 \to \RR\otimes_\kk J \to \RR_A \to \RR_B \to 0$.
Applying $\Hom_{X_A}(\KK_A,-)$ gives an exact sequence
\begin{align} \label{some_seq}
\Hom_{X_A}(\KK_A,\RR_A)
\longrightarrow
\Hom_{X_A}(\KK_A,\RR_B)
\longrightarrow
\Ext^1_{X_A}(\KK_A,\RR\otimes_\kk J).
\end{align} 
Every map $\KK_A\to\RR_B$ factors uniquely through $\KK_A\otimes_A B=\KK_B$, so the middle term is
$\Hom_{X_B}(\KK_B,\RR_B)$.
Since $\KK_A$ is flat over $A$, a locally free resolution of $\KK_A$ remains exact after reduction to $\kk$.
Consequently, $\Ext^1_{X_A}(\KK_A,\RR\otimes_\kk J)\simeq \Ext^1_X(\KK,\RR)\otimes_\kk J=0$.
Thus, the first map in \Cref{some_seq} is surjective, which implies that $i_B$ lifts to a map $i_A:\KK_A\to\RR_A$. By induction, its restriction to $\kk$ is the original inclusion $\KK\hookrightarrow\RR$.

In either case, the fiberwise injectivity criterion shows that  $i_A$ is injective with $A$-flat
cokernel; in particular, its image $\MM_A$ lies in $\Def_{z'}(A)$.
Moreover, $i_A$ is an isomorphism on $U_A$, since its cokernel its zero on the closed fiber.
Thus, on $U_A$, the map $\kappa_Ai_A^{-1}$ defines a map $\RR_A|_{U_A} \to \EE_A|_{U_A}$.
Applying $(j_{A})_\star$ extends it uniquely to a map $f_A:\RR_A\to\EE_A$.
Since $\kappa_A=f_Ai_A$ holds on $U_A$, and $\EE_A \to (j_A)_{\star}(\EE_A|_{U_A})$ is injective, the equality also holds on $X_A$.
Finally, since the restriction of $f_A$ to $\kk$ is $f_0$ we see that $\Xi_A(\MM_A,f_A) = [\EE_A\twoheadrightarrow \FF_A]$, proving surjectivity.

We are now ready to show that $\Xi: \mathrm{D} \to \Def_z$ is formally smooth. Indeed, let $A'\twoheadrightarrow A$ be a small extension, and consider the commutative diagram
\[
\begin{tikzcd}
\Spec(A) \arrow[r, "\eta"]
\arrow[d, hook] &
\mathrm{D} \arrow[d, "\Xi"] \\
\Spec(A') \arrow[r, "\xi'"'] &
\Def_z
\end{tikzcd}
\]
By surjectivity of $\Xi_{A'}$, there exists $\eta'\in \mathrm{D}(A')$ with $\Xi(\eta')=\xi'$.
Since $ \Xi: \mathrm{D}(A)/\mathrm{G}(A) \simeq \Def_z(A)$, there exists $g\in \mathrm{G}(A)$ such that $g\cdot(\eta'|_A)=\eta$.
Since $\mathrm{G}(A')\to \mathrm{G}(A)$ is surjective, we may choose a lift $g'\in \mathrm{G}(A')$ of $g$. Then $g'\cdot \eta'$ lifts $\eta$ and maps to $\xi'$. 
In particular, $\Xi$ is formally smooth.
\end{proof}

We now state two cases where the hypotheses of the theorem apply. 
The first, as mentioned earlier, is an extension of Fogarty's result of detaching hypersurfaces to Quot schemes, while the second generalizes \cite[Lemma 2.1(ii)]{Ricolfi}. 
We begin with a definition

\begin{definition} A coherent sheaf $\EE$ on a scheme $X$ is called \textbf{rigid} if $\Ext^1_X(\EE,\EE) = 0$.
\end{definition}

\begin{cor} \label{double_dual_cor}
Let $X$ be a smooth projective scheme of dimension at least two, and let $\EE$ be a locally free sheaf. Let $[\EE \twoheadrightarrow \FF]$ be a point of $\Quot_X(\EE)$ with kernel $\KK$. Then the conclusion of \Cref{double_dual_prop} holds if either
\begin{enumerate}
\item $H^1(X,\OO_X)=0$ and $\KK= \II_Z\otimes \mathscr{L}$ for a closed subscheme $Z\subseteq X$ of codimension at least two and an invertible sheaf $\mathscr{L}$, or
\item $\KK^{\vee\vee}$ is locally free, rigid, and $\codim \Supp(\KK^{\vee\vee}/\mathscr K) \geq 3$.
\end{enumerate}
\end{cor}
\begin{proof}
The first item follows from the fact that $\KK^{\vee\vee} = \mathscr{L}$ is of rank one. 
For the second item, let $n$ be the dimension of $X$. 
By Serre duality we have 
$$
\Ext^2_X(\KK^{\vee\vee}/\KK,\KK^{\vee\vee}) \cong H^{n-2}(X,(\KK^{\vee\vee})^\vee \otimes (\KK^{\vee\vee}/\KK)\otimes\omega_X)^{\vee},
$$
which vanishes by our assumptions. 
Applying $\Hom_X(-,\KK^{\vee\vee})$ to $0 \to \KK \to \KK^{\vee\vee} \to \KK^{\vee\vee}/ \KK \to 0$ gives us the exact sequence
$$
\Ext_X^1(\KK^{\vee\vee},\KK^{\vee\vee}) \to \Ext_X^1(\mathscr K,\KK^{\vee\vee}) \to \Ext_X^2(\KK^{\vee\vee}/\KK,\KK^{\vee\vee}).
$$ 
Since $\KK^{\vee\vee}$ is rigid, the leftmost group also vanishes, implying that the middle group vanishes, as required.
\end{proof}

\begin{lemma}\label{Deformation_theory_flag}
Let $X$ be a projective scheme and $\EE$ be a coherent sheaf. Consider the natural morphism
$$
\pi:
\operatorname{FlagQuot}^{P_1,P_2}_X(\EE)
\longrightarrow
\Quot^{P_1}_X(\EE).
$$
Let $z=[\EE\twoheadrightarrow\FF_1\twoheadrightarrow \FF_2]$
be a point of the flag Quot scheme and put
$\NN=\ker(\FF_1\twoheadrightarrow\FF_2)$.
If $\Hom_X(\NN,\FF_2)=\Ext^1_X(\NN,\FF_2)=0$, then $\pi$ is \'etale at $z$.
\end{lemma}
\begin{proof}
Let $\QQ=\Quot_X^{P_1}(\EE)$ and let $\MM$ be the universal quotient on $X\times \QQ$. 
By the definition of the flag Quot functor, we have
$$
\mathrm{FlagQuot}_X^{P_1,P_2}(\EE) \simeq \Quot_{X\times \QQ/\QQ}^{P_2}(\MM),
$$
under which $\pi$ is the structure morphism to $\QQ$. 
By \cite[Lemma 2.2.6]{HL}, given a small extension of local Artinian $\kk$-algebras with kernel $I$, the obstruction to lifting a deformation of $z$, relative to $\QQ$, lies in $\Ext^1_X(\NN,\FF_2)\otimes_{\kk} I$.
When such a lift exists, the set of lifts is a torsor under $\Hom_X(\NN,\FF_2)\otimes_{\kk} I$.
By assumption, both groups vanish and thus there is a unique lift.
By the infinitesimal criterion for \'etaleness, $\pi$ is \'etale at $z$.
\end{proof}

\section{Quot schemes on projective space}
We review the Hilbert polynomials of quotients of absolutely split bundles and the connectedness of their Quot schemes.
This will provide the numerical and combinatorial setup for the study of lexicographic points and their components in \Cref{sec:lex}. 
To illustrate why we focus on the absolutely split case, we give an example of a disconnected Quot scheme when the source is indecomposable of rank two.

\begin{example} \label{disconnected_ex} Let $\EE$ be the cokernel of the injective map 
$$
\phi: \OO_{\PP^2}(-2) \xrightarrow{\begin{pmatrix} x_0^3+x_1^3 & x_2^3 & x_0x_1\end{pmatrix}^T} \OO_{\PP^2}(1)^2 \oplus\OO_{\PP^2}
$$
Since the vanishing locus of the entries is empty, $\EE$ is locally free of rank two.  We claim that $\Quot_{\PP^2}^{t+1}(\EE)$ is disconnected.

Consider a point $[\OO_{\PP^2}(1)^2 \oplus \OO_{\PP^2} \twoheadrightarrow \FF]$ on the Quot scheme $\Quot_{\PP^2}^{t+1}(\OO_{\PP^2}(1)^2 \oplus \OO_{\PP^2})$ and let  $\NN \subseteq \FF$ be its maximal zero-dimensional subsheaf. 
Since the pure sheaf $\FF/\NN$ also has degree one, its support is a line $L$ in $\PP^2$. 
Since $\FF/\NN$ has rank one we must have $\FF/\NN \cong\OO_L(a)$.
Since $\OO_{\PP^2}(1)^{\oplus 2} \oplus \OO_{\PP^2}$ is globally generated, so is $\FF$, and thus $a \geq 0$. 
Comparing Hilbert polynomials we obtain $t+1 - P_{\NN} = t+1+ a$. 
In particular, $a = 0$ and $\FF \cong \OO_L$.
Thus, restriction defines a morphism
$$
\Psi:\Quot_{\PP^2}^{t+1}(\OO_{\PP^2}(1)^{\oplus2}\oplus\OO_{\PP^2}) \longrightarrow \Hilb^{t+1}(\PP^2).
$$
It has a section induced by the projection $\OO_{\PP^2}(1)^{\oplus2}\oplus\OO_{\PP^2}\twoheadrightarrow\OO_{\PP^2}$.
Moreover, $\Ext^1_{\PP^2}(\OO_{\PP^2}(1),\OO_L)=H^1(L,\OO_L(-1))=0$, so  \Cref{fibration} and the fiber description show that $\Psi$ is smooth of relative dimension zero and, since its closed fibers are points, it is an isomorphism.
Consequently, $\Quot_{\PP^2}^{t+1}(\OO_{\PP^2}(1)^{\oplus2}\oplus\OO_{\PP^2}) \simeq\Hilb^{t+1}(\PP^2)\simeq(\PP^2)^\vee.$

A quotient $\OO_{\PP^2}(1)^2 \oplus \OO_{\PP^2} \twoheadrightarrow \OO_{\PP^2} \twoheadrightarrow \OO_L$ factors through $\EE$ precisely when $x_0x_1|_L=0$.
Thus $\Quot_{\PP^2}^{t+1}(\EE)$ is the closed subscheme of $(\PP^2)^\vee$ parametrizing lines contained in $V(x_0x_1)$.
Its only points are $V(x_0)$ and $V(x_1)$; in particular, the Quot scheme is disconnected.
In fact, we have  $$\Quot_{\PP^2}^{t+1}(\EE)\cong \Spec(\kk) \sqcup\Spec(\kk).$$
\end{example}

The indecomposability of $\EE$ in the above example will follow from the fact that Quot schemes of absolutely split bundles are connected (see \Cref{Pardue_conn}). 
This connectedness can be attributed to the presence of a large group action, which allows us to degenerate to Borel-fixed points, just as in the case of Hilbert schemes. 
We discuss these ideas below, beginning with the structure of Hilbert polynomials.

\subsection{Integer partitions and connectedness}
\begin{definition}
An \textbf{integer partition} is a sequence $\lambda=(\lambda_1,\ldots,\lambda_e)$ with $\lambda_1\geq\cdots\geq\lambda_e\geq1$. 
We associate to a partition $\lambda$ the polynomial
$$
p_\lambda(t)=\sum_{i=1}^e\binom{t+\lambda_i-i}{\lambda_i-1}.
$$
The number of times a $j$ occurs in $\lambda$ is called its \textbf{multiplicity}.
We also allow the empty partition $\lambda=\varnothing$, with $p_\varnothing=0$ and $\deg(p_\varnothing)=-1$.
\end{definition}

\begin{definition} Given an integer partition $\lambda $, we define $\partial\lambda$ to be the integer partition obtained from $\lambda $ by subtracting one from each part and deleting the zero entries.

\end{definition}
\begin{remark} If $\lambda$ is an integer partition, then $\Delta p_\lambda = p_{\partial \lambda}$ where $\Delta$ is the difference operator $\Delta f(t) := f(t)-f(t-1)$.
\end{remark}

\begin{notation}\label{not:quot}
Let $n\geq 1$ and $S = \kk[x_0,\dots,x_n]$ with its standard grading. 
For a finitely generated graded $S$-module $M$, we use $\mfh_M$ to denote its Hilbert function.
Fix integers $a_1 \leq a_2 \leq \cdots \leq a_r$ and let
$$ 
E_{\bfa}=\bigoplus_{i=1}^r S(-a_i)\bfe_i, \qquad \EE_{\bfa}=\bigoplus_{i=1}^r\OO_{\PP^n}(-a_i).
$$
Here $\bfe_i$ denotes the generator of the $i$-th summand, of degree $a_i$.
For $0\leq s<r$, put $q=r-s$. 

For a partition $\lambda$ whose parts have size at most $n$, set
$$
P_{\bfa,s,\lambda}(t) = \sum_{j=q+1}^r\binom{t-a_j+n}{n}+p_\lambda(t-a_q), \qquad
\QQ_{\bfa}(s,\lambda) = \Quot^{P_{\bfa,s,\lambda}}_{\PP^n}(\EE_{\bfa}).
$$
The integer $s$ is the rank of the quotient, and $q$ is the rank of its kernel. 

We write 
$$
0 \longrightarrow \KK_\lambda \longrightarrow (\EE_\bfa)_{\QQ_{\bfa}(s,\lambda)} \longrightarrow \FF_\lambda \longrightarrow 0
$$
for the universal sequence on $\PP^n \times \QQ_{\bfa}(s,\lambda)$.
\end{notation}

We will see in \Cref{Macaulay} that every quotient of $\EE_{\bfa}$ has a Hilbert polynomial of the above form.

\begin{notation} Let $\EE_{\bfa} \twoheadrightarrow \FF$ be a quotient and $\KK$ its kernel. 
Then there is a unique saturated homogeneous submodule $K\subseteq E_{\mathbf a}$ such that $\KK = \widetilde{K}$. 
We often denote the point corresponding to $[\EE \twoheadrightarrow \FF]$ in $\Quot_{\PP^n}(\EE)$ by $[E_{\bfa}/K]$.
\end{notation}

\begin{definition} 
A \textbf{monomial} in $E_{\bfa}$ is an element $\bfx^{\bfu}\bfe_i$, where $\bfx^{\bfu}$ is a monomial in $S$; it has degree $|\bfu| + a_i$. 
A submodule $K$ is said to be \textbf{monomial} if it is generated by monomials. In particular, we have a decomposition $K = \bigoplus_{i=1}^r I_i(-a_i)\bfe_i$ where each $I_i \subseteq S$ is a monomial ideal.
\end{definition}

\begin{definition} 
We put the standard \textbf{lexicographic} order on $S$, namely $\bfx^{\bfu} > \bfx^{\bfv}$ if the first non-zero entry of $\bfu-\bfv$ is positive. 
We also assume $E_{\bfa}$ has the \textbf{position-over-term} order, namely $\bfx^{\bfu}\bfe_i > \bfx^{\bfv}\bfe_j$ if $i < j$ or  ($i =j$ and $\bfx^{\bfu} > \bfx^{\bfv}$). 
A \textbf{lexicographic subspace} of $(E_{\mathbf a})_d$ is spanned by an initial segment of its monomials in decreasing order. 
A homogeneous submodule $K\subseteq E_{\bfa}$ is \textbf{lexicographic} if $K_d$ is a lexicographic subspace for every $d$. 
\end{definition}

\begin{definition} 
Let $\GL(m)$ denote the group of invertible $m\times m$ matrices and $B(m)$ the subgroup of upper-triangular matrices. 
Similarly, let $\GL(E_{\mathbf a})$ denote the group of graded $S$-module automorphisms of $E_{\mathbf a}$ and $\mathrm{B}(E_{\bfa})$ the subgroup of graded upper-triangular automorphisms. 

There is a natural action of $\mathrm{G}_{\bfa} := \GL(n+1) \rtimes \GL(E_{\bfa})$ on $E_{\bfa}$ given by $(g,A)\cdot\left(\sum_i h_i \bfe_i \right) = A\left(\sum_i g(h_i)\bfe_i \right)$.
The \textbf{Borel subgroup} is $\mathrm{B}_{\bfa}=\mathrm{B}(n+1)\rtimes\mathrm{B}(E_{\bfa})$. Thus, $\mathrm{G}_{\bfa}$ and $\mathrm{B}_{\bfa}$ act on $\QQ_{\mathbf a}(s,\lambda)$. 
A submodule $K\subseteq E_{\bfa}$ or its corresponding point on the Quot scheme is said to be \textbf{Borel-fixed} if it is invariant under $\mathrm{B}_{\bfa}$.
\end{definition}

We recall the following classification of saturated Borel-fixed submodules.

\begin{prop}[{\cite[Proposition 6]{Pardue}}] \label{Pardue_Borel}
\label{Borel_char}
A saturated homogeneous submodule $K\subseteq E_{\bfa}$ is Borel-fixed if and only if the following conditions hold:
\begin{enumerate}
\item $K=\bigoplus_{i=1}^r I_i(-a_i)\bfe_i$, where each $I_i$
is a monomial ideal;
\item each $I_i$ is Borel-fixed;
\item $I_1\supseteq I_2 \supseteq \cdots \supseteq I_r$.
\end{enumerate}
\end{prop}

Since $B_{\mathbf a}$ is connected and solvable, the Borel fixed point theorem gives

\begin{cor}\label{cor:borel-intersections}
Every irreducible component of $\QQ_{\bfa}(s,\lambda)$, and every nonempty intersection of irreducible components, contains a Borel-fixed point.
\end{cor}

Macaulay's theorem for modules asserts that every homogeneous submodule of $E_{\mathbf a}$ has a unique lexicographic submodule with the same Hilbert function. 
If we restrict to saturated submodules, it turns out that the Hilbert polynomial is precisely of the form $P_{\bfa,s,\lambda}$.

 \begin{prop}\label{Macaulay}
Let $P \in \bbQ[t]$ be a numerical polynomial.
The following statements are equivalent:
\begin{enumerate}
\item There is a quotient $\EE_{\bfa}\twoheadrightarrow\FF$ with nonzero kernel and Hilbert polynomial $P$.
\item There are an integer $0\leq s<r$ and a partition $\lambda$, whose parts have size at most $n$, such that $P=P_{\bfa,s,\lambda}$.
\end{enumerate}
When these conditions hold, the pair $(s,\lambda)$ is unique and $s=\mathrm{rank}(\FF)$.
\end{prop}

\begin{proof}
Assume that (i) is true.
By Macaulay's theorem for modules \cite[Proposition 2]{Pardue}, there is a saturated lexicographic submodule $K\subseteq E_{\bfa}$ such that $E_{\bfa}/K$ has Hilbert polynomial $P$.
Let $K=\bigoplus_{i=1}^r I_i(-a_i)\bfe_i$ and let $q$ be the largest index for which $I_q\neq0$.
Since $I_q$ contains non-zero monomials in sufficiently high degree, the position-over-term lexicographicity implies that $K$ contains every monomial in that degree in all of the preceding summands.
In particular, $I_i$ contains all monomials of sufficiently large degree for $i<q$, which, since $I_i$ is saturated, gives $I_i=S$.
The ideal $I_q$ is nonzero, saturated, and lexicographic. 
Thus, Macaulay's theorem implies that $S/I_q$ has Hilbert polynomial $p_\lambda$ for a unique partition $\lambda$ whose parts have size at most $n$ (see \cite{SS23}). 
Hence
$$
E_{\bfa}/K \simeq (S/I_q)(-a_q)\bfe_q \oplus \bigoplus_{j=q+1}^r S(-a_j)\bfe_j.
$$
Taking Hilbert polynomials gives $P=P_{\bfa,r-q,\lambda}$, proving (ii).

Conversely, suppose that (ii) holds and put $q=r-s$.
By Macaulay's theorem, there is a closed subscheme $Z \subseteq \PP^n$ with Hilbert polynomial $p_\lambda$.
Consider the projection
$$
\EE_{\bfa}\twoheadrightarrow \OO_Z(-a_q)\oplus \bigoplus_{j=q+1}^r\OO_{\PP^n}(-a_j).
$$
Its kernel is nonzero and its Hilbert polynomial
is $P_{\bfa,s,\lambda}$, proving (i).
Finally, the leading coefficient of $P$ determines $s$, and uniqueness of the Macaulay representation then determines $\lambda$.
\end{proof}

\begin{definition} \label{lexmodule}

Let $\lambda$ be a partition with parts of size at most $n$, and let $m_j$ denote the multiplicity of the part $j$ in $\lambda$.
The saturated \textbf{lexicographic ideal} $L(\lambda)\subseteq S$ is
$$
L(\lambda)=\left(
x_0^{m_n+1},\,
x_0^{m_n}x_1^{m_{n-1}+1},\, \ldots,\,
x_0^{m_n}x_1^{m_{n-1}} \cdots
x_{n-3}^{m_3}x_{n-2}^{m_2+1},\,
x_0^{m_n}x_1^{m_{n-1}} \cdots x_{n-2}^{m_2}x_{n-1}^{m_1}
\right).
$$
The quotient $S/L(\lambda)$ has Hilbert polynomial $p_\lambda$ and defines a point in $\Hilb^{p_\lambda}(\PP^n)$.

For $0\leq s<r$, let $q=r-s$. The \textbf{lexicographic submodule} corresponding to $(s,\lambda)$ is
$$
L_{\bfa}(s,\lambda) = \bigoplus_{i<q}S(-a_i)\bfe_i \oplus L(\lambda)(-a_q)\bfe_q \oplus\bigoplus_{j>q}0\bfe_j.
$$
By the proof of \Cref{Macaulay}, this is the unique saturated lexicographic submodule of $E_{\bfa}$
whose quotient has Hilbert polynomial $P_{\bfa,s,\lambda}$.
The corresponding point of $\QQ_{\bfa}(s,\lambda)$ is called the \textbf{lexicographic point}.
When this point is smooth, its unique irreducible component is called the \textbf{lexicographic component} and is denoted by $\LL_{\bfa}(s,\lambda)$.
\end{definition}

\begin{remark} \label{mingen} Let $\ell(\lambda)$ denote the \textbf{length} of $\lambda$. 
If $\lambda \ne \varnothing$, the ideal $L(\lambda)$ has exactly $n-\lambda_{\ell(\lambda)}+1$ minimal generators.
\end{remark}

Building on Reeves' proof of the connectedness of $\Hilb^{p_\lambda}(\PP^n)$, Pardue showed that

\begin{thm}[{\cite{Pardue}}] \label{Pardue_conn} The Quot schemes $\QQ_{\bfa}(s,\lambda)$ are rationally chain connected.
\end{thm}

\begin{definition} \label{radius_defn}
The \textbf{incidence graph} of $\QQ_{\bfa}(s,\lambda)$ has the irreducible components of $\QQ$ as vertices, with an edge between distinct components whenever they intersect.
Write $d(C,D)$ for the length of a shortest path between components $C$ and $D$.
For closed points $z,z' \in\QQ_{\bfa}(s,\lambda)$, define
$$
d(z,z')=
\min \{d(C,D): z \in C, \, z' \in D\},
$$
where $C$ and $D$ range over the irreducible components of $\QQ_{\bfa}(s,\lambda)$.
The \textbf{lexicographic distance} of $z$ is $d_{\text{lex}}(z) := d(z,L_{\bfa}(s,\lambda))$.
The \textbf{eccentricity} of a component $D$ is
$r_D=\max_C d(C,D)$, and the \textbf{radius} of $\QQ_{\bfa}(s,\lambda)$ is
$$
\rad(\QQ_{\bfa}(s,\lambda))=\min_D r_D,
$$
where both extrema range over the irreducible components of $\QQ_{\bfa}(s,\lambda)$. 
\end{definition}

\begin{thm}[{\cite{Reeves_radius, Pardue}}] \label{Reeves_Pardue} 
The lexicographic distance of any point on the Quot scheme $\QQ_{\bfa}(s,\lambda)$ is at most $\deg(P_{\bfa,s,\lambda}) +1$. 
\end{thm}

Pardue proved this by constructing curves on the Quot scheme using polarization and specialization. 
We will need a special case of this construction, called a distraction.

\begin{definition}\label{def:distraction}
Let $c_{ij}\in\kk$ be a collection of scalars called \textbf{distraction parameters}.
For a monomial $\bfx^{\boldsymbol{\alpha}} \in S$, define
$$
D_t(\bfx^{\boldsymbol{\alpha}}) = x_n^{\alpha_n}\prod_{i=0}^{n-1}\prod_{j=0}^{\alpha_i-1}(x_i-tc_{ij}x_n) \in S[t].
$$
If $I=(\boldsymbol{m}_1,\dots,\boldsymbol{m}_p)$ is a monomial ideal, define $D_t(I)=(D_t(\boldsymbol{m}_1),\ldots,D_t(\boldsymbol{m}_p))$. 
The ideal $D_t(I)|_{t=1} \subseteq S$ is called a \textbf{distraction} of $I$.
For a monomial submodule
$K=\bigoplus_{i=1}^r I_i(-a_i)\bfe_i\subseteq E_{\bfa}$,
define
$$
D_t(K)=\bigoplus_{i=1}^r D_t(I_i)(-a_i)\bfe_i \subseteq E_{\bfa}\otimes_\kk \kk[t].
$$
We call $D_t(K)|_{t=1}$ a \textbf{distraction} of $K$.
\end{definition}

\begin{lemma}\label{lemma_distraction}
Let $[E_{\bfa}/K]$ be a monomial point of $\QQ_{\bfa}(s,\lambda)$. Then we have a morphism
$$
\gamma_K:\mathbf{A}^1_\kk\longrightarrow\QQ_{\bfa}(s,\lambda),
\qquad t \mapsto [E_{\bfa}/D_t(K)].
$$
If $K$ is Borel-fixed, the distraction parameters may be chosen so that all the non-zero fibers of $D_t(K)$ are direct sums of radical ideals. The parameters may be chosen independently in different summands.
Moreover, using the same parameters preserves inclusions, namely,  if $K\subseteq K'$ are monomial submodules of $E_{\bfa}$, then $D_t(K)\subseteq D_t(K')$.
\end{lemma}

\begin{proof} Let $K=\bigoplus_{i=1}^r I_i(-a_i)\bfe_i$ and note that  distractions are defined summandwise. Thus, it suffices to prove the theorem for monomial ideals in $S$. 
Thus, the lemma follows from  \cite[Proposition~2.9, Corollary~2.10, Proposition~2.26]{BCR05}.
\end{proof}

\section{Lexicographic points and components}\label{sec:lex}
We determine when the lexicographic point is smooth. 
For $n\geq2$, we describe the lexicographic component by a single relative Quot construction over a space of residual flags. We begin by recalling the construction of residual flags from \cite{SS23}.

\subsection{Residual flags and rank-one kernels}

\begin{definition}[{\cite[Definition~1.4]{SS23}}] \label{residual_flag}
Let $\mu=(d_1^{c_1},\ldots,d_e^{c_e})$, where $n\geq d_1>\cdots>d_e\geq1$.
For a $\kk$-scheme $T$, a \textbf{residual flag of type $\mu$ over $T$} is a chain of closed subschemes
$$
\varnothing = Y_{e+1}\subseteq Y_e\subseteq\cdots\subseteq Y_1
\subseteq\PP^n_T
$$
such that each $Y_i$ is flat over $T$ and there exist a relative $d_i$-plane $\Lambda_i\subseteq\PP^n_T$ containing $Y_i$ and a relative effective Cartier divisor $D_i\subseteq\Lambda_i$ of degree $c_i$ satisfying
$$
\II_{Y_i/\Lambda_i}
=\II_{D_i/\Lambda_i}\II_{Y_{i+1}/\Lambda_i}.
$$
The functor $T \mapsto \{\text{residual flags of type $\mu$ in $\PP^n_T$}\}$ is represented by the scheme $\Flag(\mu;\PP^n)$. By abuse of terminology, we call $Y_1$ the residual flag.
\end{definition}

\begin{remark} \label{universal_flag}
By \cite[Proposition~1.11]{SS23}, the scheme $\Flag(\mu;\PP^n)$ is smooth, irreducible, and projective,
and is built from projective bundles over partial flag varieties. It carries a universal residual flag
$$
\varnothing=\mathcal Y_{e+1}\subseteq\mathcal Y_e
\subseteq\cdots\subseteq\mathcal Y_1
\subseteq\PP^n\times\Flag(\mu;\PP^n).
$$
\end{remark}

\begin{thm}[{\cite[Theorems~1.4 and~4.1]{RS97}}]\label{hilbert_lex}
Write $\lambda=(\mu,1^b)$, where $\mu$ has no parts equal to one. 
The lexicographic point is smooth on $\Hilb^{p_\lambda}(\PP^n)$, and its component is the closure
of the locus of subschemes $Y\sqcup\{p_1,\ldots,p_b\}$, where $Y$ is a residual flag of type $\mu$ and $p_1,\ldots,p_b$ are distinct points outside $Y$.
\end{thm}

\begin{lemma}\label{negative_cohomology}
Let $n\geq 2$, and let $Y$ be a residual flag whose partition has no parts equal to one. 
Then $H^0(Y,\OO_Y(v))=0$ for $v<0$ and $H^0(Y,\OO_Y)=\kk$. Consequently, $H^1(\PP^n, \II_Y(v))=0$ for all $v$. 
\end{lemma}
\begin{proof}
For a closed subscheme $Z\subseteq\PP^n$, let $R_Z=S/I_Z$, where $I_Z$ is its saturated homogeneous ideal.
The residual construction gives exact sequences
$$
0\longrightarrow R_{Y_{i+1}}(-c_i) \longrightarrow R_{Y_i}\longrightarrow R_{D_i} \longrightarrow 0.
$$
By our assumptions, each $R_{D_i}$ is Cohen-Macaulay of dimension $d_i \geq 2$. 
Starting with $R_{Y_e} = 0$ and inducting downwards on $i$, the exact sequence implies that the depth of $R_Y$ is at least $2$.
Consequently, $H^0_{\mfm}(R_Y)= H^1_{\mfm}(R_Y)=0$.
Thus, $R_Y\simeq\bigoplus_{v\in\mathbf Z}H^0(Y,\OO_Y(v))$ which implies $H^0(\OO_Y(v))=0$ for all $v < 0$ and equals $\kk$ for $v=0$.
We also have $H^1(\II_Y(v))\simeq H^1_{\mfm}(R_Y)_v=0$ for all $v$. 
\end{proof}

\begin{prop} \label{P1_trivial}
Let $\lambda=(n^c,\nu)$ and assume $\nu$ has no part equal to $n$.
There is a natural isomorphism
$$
\Hilb^{p_\nu}(\PP^n)\times \PP\left(\bigoplus_{j=1}^r S_{c+a_1-a_j}\right) \xrightarrow{\ \sim\ } \QQ_{\bfa}(r-1,\lambda),
\qquad
(W,[f]) \mapsto \left[\EE_{\bfa}\twoheadrightarrow \EE_{\bfa}/f(\II_W(-a_1-c))\right].
$$
In particular, $\QQ_{\bfa}(r-1,\lambda)$ is smooth if and only if $\Hilb^{p_\nu}(\PP^n)$ is smooth.
\end{prop}

\begin{proof} 
Let $d=a_1+c$, let $V=\Hom_{\PP^n}(\OO_{\PP^n}(-d),\EE_{\bfa})$ and let $p_T:\PP^n_T\to T$ be the projection.
A $T$-point of $\Hilb^{p_\nu}(\PP^n)\times\PP(V)$ consists of a flat family $W\subseteq\PP^n_T$ with Hilbert polynomial $p_\nu$ and a line subbundle $\LLL\hookrightarrow V\otimes_\kk\OO_T$.
The inclusion $\mathscr L\hookrightarrow V\otimes_\kk\OO_T$ and the canonical evaluation map $\OO_{\PP^n}(-d)\otimes_\kk V\rightarrow\EE_{\bfa}$ induces a morphism $f:\OO_{\PP^n_T}(-d)\otimes p_T^\star\LLL \rightarrow(\EE_{\bfa})_T$.
Define the map on $T$-points by
$$
(W,\LLL\hookrightarrow V\otimes_\kk\OO_T) \mapsto \left[(\EE_{\bfa})_T\twoheadrightarrow \frac{(\EE_{\bfa})_T}{f(\II_W(-d)\otimes p_T^\star\LLL)}
\right].
$$
On every geometric fiber, $f$ is a nonzero map from a line bundle to a vector bundle on $\PP^n$, hence is injective.
Its restriction to $\II_W(-d)\otimes p_T^\star\LLL$ is therefore also fiberwise injective. 
Since both source and target are flat over $T$, the fiberwise flatness criterion shows that its cokernel is flat over $T$.
The Hilbert polynomial of each fiber is
$$
P_{\EE_{\bfa}}(t)-\binom{t-d+n}{n}+p_\nu(t-d) =
P_{\EE_{\bfa}}(t) + p_\lambda(t-a_1)-\binom{t-a_1+n}{n} = P_{\bfa,r-1,\lambda}(t).
$$
The construction commutes with arbitrary base change, so these maps on $T$-points define the required morphism.

Consider a closed point $z=[\EE_{\bfa}\twoheadrightarrow\FF]$with kernel $\KK$. 
Since $\KK$ is torsion-free of rank one, its double dual is a line bundle, say $\KK^{\vee\vee} \simeq \OO_{\PP^n}(-e)$, and $\KK^{\vee\vee}/\KK$ is supported in codimension at least two.
On the other hand,
$$
P_{\KK}(t) = P_{\EE_{\bfa}}(t)-P_{\bfa,r-1,\lambda}(t) = \binom{t-d+n}{n}-p_\nu(t-d).
$$
Since $\deg p_\nu\leq n-2$, comparison of the coefficients of $t^{n-1}$ gives $e=d$.
The inclusion $\KK\hookrightarrow\KK^{\vee\vee}$ therefore identifies $\KK$ with $\II_W(-d)$ for a unique subscheme $W\subseteq\PP^n$ of codimension at least two.
The exact sequence $0 \to \KK\to\OO_{\PP^n}(-d) \to\OO_W(-d)\to 0$ then gives $P_W=p_\nu$.
Taking double duals of $\KK\hookrightarrow\EE_{\bfa}$ gives an inclusion $\OO_{\PP^n}(-d) \hookrightarrow \EE_{\bfa}$.
Changing the identification of $\KK^{\vee\vee}$ with $\OO_{\PP^n}(-d)$ rescales the polynomial $f$.
Thus, the point $z$ determines a unique pair $(W,[f])$, proving bijectivity on geometric points.

Since the morphism is proper and bijective on geometric points, it remains to prove it is \'etale.
This follows from arguing as in the proof of \Cref{double_dual_prop} (i) with $\RR=\OO_{\PP^n}(-d)$.
Indeed, for a local $\kk$-Artinian algebra $A$, the elements of $\mathrm{D}(A)$ are pairs $(\II_{W_A}(-d),f_A)$ deforming $(\II_W(-d),f)$. 
The group $\mathrm G(A)=1+\mathfrak m_A$ fixes $\II_{W_A}(-d)$ and rescales $f_A$.
Thus $\mathrm D/\mathrm G$ is precisely the local deformation functor of $(W,[f])$ in the product.
The isomorphism $\mathrm D/\mathrm G\simeq\Def_z$ proved there is induced by our morphism. 
Thus, the morphism is \'etale, as required.
\end{proof}

\subsection{Smoothness of the lexicographic point}
Throughout this subsection, write $\lambda=(n^c,\nu)$, where $\nu$ has no part equal to $n$, and put $d=a_q+c$. 
Decompose $\EE_{\bfa}=\EE_1\oplus\EE_2$, with $\EE_1=\bigoplus_{i<q}\OO_{\PP^n}(-a_i)$ and $\EE_2=\bigoplus_{j=q}^r\OO_{\PP^n}(-a_j)$, and let $\QQ_2=\Quot_{\PP^n}^{P_{\bfa,s,\lambda}}(\EE_2)$.
By \Cref{P1_trivial}, we have an isomorphism
$$
\QQ_2\cong \Hilb^{p_\nu}(\PP^n)\times \PP\left(\bigoplus_{j=q}^r S_{d-a_j}\right).
$$
Under this identification, a point $(W,[\bff])$ corresponds to the quotient in the exact sequence
\begin{align} \label{quot_prod}
0\longrightarrow\II_W(-d) \xrightarrow{\bff}\EE_2 \longrightarrow\FF\longrightarrow0.
\end{align}
Let $Z_\nu$ be the scheme defined by $L(\nu)$, and let $z=([Z_\nu],[(x_0^c,0,\ldots,0)])$.
Since $L(\lambda)=x_0^cL(\nu)$, $z$ is the lexicographic point of $\QQ_2$.
Moreover, $\QQ_2$ is smooth at $z$. 
Indeed, for $n\geq2$ this follows from \Cref{hilbert_lex}, while for $n=1$ we have $\nu=\varnothing$ and $\QQ_2$ is a projective space.

Finally, by \Cref{restriction_open}, restriction induces a morphism $\Psi:\UU_{\EE_2}\to\QQ_2$.
Projection onto $\EE_2$ gives a section $\sigma$ taking $z$ to the lexicographic point of $\QQ_{\bfa}(s,\lambda)$. 
Thus, the above discussion and \Cref{fibration}(ii) reduce smoothness at $\sigma(z)$ to local constancy near $z$ of the function
\begin{align}\label{function_const}
[\EE_2\twoheadrightarrow\FF] \longmapsto h^0(\EE_1^\vee\otimes\FF) =\sum_{i<q}h^0(\FF(a_i)).
\end{align}

\begin{prop}\label{lex_smooth}
If $n\geq2$, the lexicographic point of $\QQ_{\bfa}(s,\lambda)$ is smooth.
\end{prop}

\begin{proof}
By the preceding reduction, it suffices to prove that \eqref{function_const} is locally constant near $z$.
For $i<q$, taking cohomology of \Cref{quot_prod} and of the ideal sheaf sequence of $W$,
we obtain
\begin{align*}
h^0(\FF(a_i)) &= h^0(\EE_2(a_i)) - h^0(\II_W(a_i-d)) + h^1(\II_W(a_i-d)) \\
& = h^0(\EE_2(a_i)) -h^0(\OO_{\PP^n}(a_i-d)) +h^0(\OO_W(a_i-d))
\end{align*}
for all $i < q$.
Thus, to establish the constancy of \Cref{function_const}, it is enough to establish, for each $i < q$, the constancy of $h^0(\OO_W(a_i-d))$ for $W$ in a neighbourhood of $[L(\nu)]$.

Let $\nu=(\mu,1^b)$ with $\mu$ having no parts equal to one.
The inclusion $L(\nu)\subseteq L(\mu)$ gives
$$
0\longrightarrow\widetilde{L(\mu)/L(\nu)} \longrightarrow\OO_{Z_\nu}\longrightarrow\OO_{Z_\mu} \longrightarrow 0,
$$
with the first term having length $b$. 
By \Cref{negative_cohomology}, we have
$$
h^0(\OO_{Z_\nu}(a_i-d))=
\begin{cases}
b, & a_i-d <0 \text{ or } \mu = \varnothing,\\
b+1, & a_i-d =0 \text{ and } \mu \ne \varnothing.
\end{cases}
$$
By \Cref{hilbert_lex}, these are also the values for a general point of the lexicographic component $\Hilb^{p_\nu}(\PP^n)$. 
It follows from upper semicontinuity that $h^0(\OO_W(a_i-d))$ is constant for all $W$ in a neighbourhood of $[L(\nu)]$, as required.
\end{proof}

On the other hand, when $n=1$ the lexicographic point can be singular and can lie on the intersection of two components.

\begin{example} \label{non_distinguish} 
Consider the Quot scheme 
$$
\QQ = \Quot_{\PP^1}^{t+2}
\left(
\OO_{\PP^1}(2) \oplus \OO_{\PP^1}(1)\oplus\OO_{\PP^1}
\right).
$$
We will show that this Quot scheme has exactly two components and both of them contain the unique Borel-fixed point. Consider the exact sequence associated to a closed point
$$
0 \longrightarrow \KK \longrightarrow \OO_{\PP^1}(2) \oplus \OO_{\PP^1}(1)\oplus\OO_{\PP^1} 
\longrightarrow \FF \longrightarrow 0.
$$
The Hilbert polynomial of $\KK$ is $(3t+6) - (t+2) = 2t+4$. 
Since $\KK$ is locally free of rank two, the only possibilities are $\KK \cong \OO_{\PP^1}(1)^{\oplus 2}$ or $\KK \cong \OO_{\PP^1}(2) \oplus \OO_{\PP^1}$.

Consider the two Quot schemes
$\QQ_1=\Quot_{\PP^1}^{t+2}(\OO_{\PP^1}(1) \oplus \OO_{\PP^1})$ and $\QQ_2=\Quot_{\PP^1}^{1}(\OO_{\PP^1}(2) \oplus \OO_{\PP^1}(1))$.
We have natural closed immersions $\iota_1:\QQ_1 \to \QQ$ given by 
$$
[\OO_{\PP^1}(1) \oplus \OO_{\PP^1} \overset{f}{\twoheadrightarrow} \FF] \mapsto 
[\OO_{\PP^1}(2) \oplus \OO_{\PP^1}(1) \oplus \OO_{\PP^1} \overset{\boldsymbol{0} \, \oplus f}{\twoheadrightarrow}\FF]
$$
and  $\iota_2:\QQ_2 \to \QQ$ given by
$$
[\OO_{\PP^1}(2)\oplus\OO_{\PP^1}(1) \overset{f}{\twoheadrightarrow}  \FF] \mapsto
[\OO_{\PP^1}(2)\oplus\OO_{\PP^1}(1)\oplus\OO_{\PP^1} \overset{f \oplus \, \text{id}}{\twoheadrightarrow} \FF \oplus \OO_{\PP^1}].
$$
If $\KK\simeq\OO_{\PP^1}(2) \oplus\OO_{\PP^1}$, the summand $\OO_{\PP^1}(2)$ must map isomorphically to the first summand of $\EE_\bfa$. 
In particular, $\EE_{\bfa}/\KK$ belongs to $\iota_1(\QQ_1)$.
On the other hand, if $\KK\simeq\OO_{\PP^1}(1)^{\oplus2}$, then $\Hom_{\PP^1}(\KK,\OO_{\PP^1})=0$ and thus $[\EE_\bfa/\KK]$ must belong to $\iota_2(\QQ_2)$.
In particular, set-theoretically, we have $\QQ=\iota_1(\QQ_1) \cup \iota_2(\QQ_2)$.

Note that the kernel of a quotient in $\QQ_1$ must have Hilbert polynomial $t+1$. 
The only option is $\OO_{\PP^1}$ and thus
$$
\QQ_1 \cong \PP(\Hom_{\PP^1}(\OO_{\PP^1},\OO_{\PP^1}(1) \oplus \OO_{\PP^1}) = 
\PP(H^0(\PP^1,\OO_{\PP^1}(1) \oplus \OO_{\PP^1})) \simeq \PP^2
$$
is irreducible.
On the other hand $\QQ_2 \simeq \PP(\OO_{\PP^1}(-2)\oplus \OO_{\PP^1}(-1))$ is also irreducible of dimension two. 
Thus, $\QQ$ is a union of two irreducible components.
Finally, the intersection consists of quotients of the form
$$
\OO_{\PP^1}(2)\oplus\OO_{\PP^1}(1)\oplus\OO_{\PP^1} \xrightarrow{\mathbf{0} \, \oplus \, \mathrm{ev}_p \, \oplus \, \mathrm{id}} \kk_p \oplus \OO_{\PP^1} 
$$
with $p \in \PP^1$. 
In particular, the intersection is parameterized by a $\PP^1$ and contains the unique Borel-fixed point $\kk_{[0:1]} \oplus \OO_{\PP^1}$.
\end{example}

Nonetheless, we still obtain the desired radius bound.

\begin{cor} \label{radius_bound} The radius of the Quot scheme $\QQ_{\bfa}(s,\lambda)$ is at most $\deg(P_{\bfa,s,\lambda}) +1$. 
\end{cor}
\begin{proof}
Suppose first that $n\geq2$. By \Cref{lex_smooth}, the lexicographic point is smooth and hence belongs
to a unique component $\LL$. By \Cref{Reeves_Pardue}, it follows that every component has distance at most $\deg(P_{\bfa,s,\lambda})+1$ from $\LL$, giving the bound.

Suppose $n=1$, so $\lambda=(1^b)$.
The natural inclusion of the first $q$ summands of $\EE_{\bfa}$ induces a morphism
$$
\Quot_{\PP^1}^{\,b} \left(\bigoplus_{i=1}^q\OO_{\PP^1}(-a_i)\right) \longrightarrow \QQ_{\bfa}(s,\lambda),
$$
obtained by adjoining the identity quotient on the remaining summands. 
By \cite[Corollary 4.7]{BFP20}, the source is irreducible and by \Cref{Borel_char} the image contains every Borel-fixed point. 
Let $\HH$ be any irreducible component containing the image. 
By \Cref{cor:borel-intersections}, every irreducible component must meet $\HH$; in particular, the radius is at most $1$.
\end{proof}

\begin{thm}\label{lex_smooth_all}
The lexicographic point of $\QQ_{\bfa}(s,\lambda)$ is singular
if and only if $n=1$, $\lambda=(1^b)$ with $b>0$, $q\geq2$,
and there is an index $j>q$ such that
$$
a_1+2\leq a_j\leq a_q+b.
$$
\end{thm}

\begin{proof}
By \Cref{lex_smooth}, we only need to check the case of $n=1$.
Let $\bff=(f_q,\ldots,f_r)$ and consider the associated quotient $0\to\OO_{\PP^1}(-d)\xrightarrow{\bff}\EE_2\to\FF_{\bff}\to0$.
For $i<q$, the long exact sequence in cohomology and Serre duality gives 
\begin{align} \label{curve_h0}
h^0(\FF_{\bff}(a_i)) =h^0(\EE_2(a_i))-h^0(\OO_{\PP^1}(a_i-d)) +h^1(\OO_{\PP^1}(a_i-d))-\mathrm{rank}\,\theta_i(\bff),
\end{align}
where   
$$
\theta_i(\bff):\bigoplus_{j=q}^r S_{a_j-a_i-2} \longrightarrow S_{d-a_i-2},\qquad (h_j)\longmapsto\sum_{j=q}^r f_jh_j.
$$
The first three terms of \Cref{curve_h0} are independent of $\bff$. Thus, it suffices to determine when
\begin{align} \label{constant_fun} 
\bff \mapsto \bigoplus_{i<q}\theta_i(\bff)
\end{align}
has constant rank in a neighbourhood of $z$.

If $q=1$, \Cref{constant_fun} is an empty sum. 
If $b=0$, then $z = [1,0,\dots,0]$ and any $[\bff]$ near $z$ will have $f_q$ being a non-zero constant; in particular, every $\theta_i(\bff)$ is surjective. 
In either case, the rank of $\theta_i(\bff)$ is constant as required.

Suppose $q\geq2$ and $b>0$. Assume that there is no $j$ satisfying the hypothesis of the theorem. 
The map $\theta_i(\bff)$ can be non-zero on its $j$-th component only if $a_i +2 \leq a_j \leq a_q+b$. 
Since $a_i +2 \geq a_1+2$, our assumptions imply that this cannot occur for $j>q$. 
In particular, the image of $\theta_i(\bff)$ is $f_qS_{a_q-a_i-2}$ which has constant rank in a neighbourhood of $[(x_0^b,0,\dots,0)]$.

Conversely, assume there exists $j>q$ with $a_1+2 \leq a_j \leq d$.
Consider $\bff(t)$ with entries $f_j= tx_1^{d-a_j}$, $f_q=x_0^b$ and all other $f_i = 0$.
Note that for all $i$, the image of $\theta_i(\bff(0))$, which is $x_0^b S_{a_q-a_i-2}$, is always contained in the image of $\theta_i(\bff(t))$.
However, for $t \ne 0$, the image of $\theta_1(\bff(t))$ also contains $x_1^{d-a_1-2}$.
In particular, $\sum_{i<q}\mathrm{rank}(\theta_i(\bff))$ is larger for $t\ne 0$ than for $t=0$, as required.
\end{proof}

\subsection{Residual quotients and the lexicographic component}

For the remainder of this section, assume $n\geq2$ and write $\lambda=(\mu,1^b)$, where $\mu$
has no parts equal to one.

\begin{definition}
We distinguish the occurrences of $a_q$ on either side of the index $q$ by putting $q_0=\#\{i\leq q:a_i=a_q\}$ and $s_0=\#\{j>q:a_j=a_q\}$.
Thus $a_q$ has multiplicity $q_0+s_0$ in $\bfa$, with $q_0\geq1$ and $s_0\geq0$. 
Let $\mathbf{G} :=\Gr(s_0+1,q_0+s_0)$ be the Grassmannian parameterizing rank $(s_0+1)$ quotients of $\kk^{q_0+s_0}$.
\end{definition}

\begin{prop}\label{core_isom}
Let $\mu=(n^c,\nu)$ with $\nu$ having no parts of size $n$. 
Consider the locally free sheaf
$$
\mathscr{V}_c=(S_c\otimes_\kk\mathscr{Q})\oplus \bigoplus_{a_j>a_q}S_{a_q+c-a_j}\otimes_\kk\OO_{\mathbf{G}}.
$$
on $\mathbf{G}$ where  $\mathscr{Q}$ denotes the universal quotient.
If $\mu$ is non-empty, there is a natural isomorphism
$$
\QQ_{\bfa}(s,\mu) \simeq \Flag(\nu;\PP^n)\times\PP_{\mathbf{G}}(\mathscr{V}_c).
$$
For $\mu=\varnothing$, one has $\QQ_{\bfa}(s,\varnothing)  \simeq \Gr(s_0,q_0+s_0)$.
In particular, $\QQ_{\bfa}(s,\mu) $ is smooth and irreducible.
\end{prop}
\begin{proof}
Suppose first that $\mu\ne\varnothing$, and put $d=a_q+c$ and $I=\{i:a_i=a_q\}$.
For each subset $J\subseteq I$ of cardinality $s_0+1$, write $\EE_{\bfa}=\EE'_J\oplus\EE_J$, where
$\EE_J= \bigoplus_{j\in J}\OO_{\PP^n}(-a_q) \oplus\bigoplus_{a_j>a_q}\OO_{\PP^n}(-a_j)$.
Let $\UU_{J} \subseteq\QQ_{\bfa}(s,\mu)$ be the open locus of quotients whose restriction to $\EE_J$ is surjective, and let $ \QQ_J=\Quot_{\PP^n}^{P_{\bfa,s,\mu}}(\EE_J)$ (see \Cref{restriction_open}).
Since $\EE_J$ has rank $s+1$, \Cref{P1_trivial} and \cite[Theorem~3.2]{SS23} give
$$
\QQ_J\simeq \Flag(\nu;\PP^n) \times
\PP\left( \bigoplus_{j\in J}S_c \oplus\bigoplus_{a_j>a_q}S_{d-a_j}\right).
$$
For a quotient $[\EE_J\twoheadrightarrow\FF] \in \QQ_J$, its kernel $\KK$ is torsion-free of rank one, with Hilbert polynomial $P_{\KK}(t)=\binom{t-d+n}{n}-p_\nu(t-d)$.
As in \Cref{P1_trivial}, this gives $\KK^{\vee\vee}\simeq\OO_{\PP^n}(-d)$ and $\KK\simeq\II_Y(-d)$, where $Y$ is the largest member of the corresponding residual flag.
By \Cref{negative_cohomology}, for every summand $\OO_{\PP^n}(-a_i)$ of $\EE_J'$ we have
$H^0(\II_Y(a_i-d))=H^1(\II_Y(a_i-d))=0$.
Taking the direct sum over the summands of $\EE'_J$ gives
$$
\Hom_{\PP^n}(\EE'_J,\II_Y(-d)) = \Ext^1_{\PP^n}(\EE'_J,\II_Y(-d)) =0.
$$
Thus every map $\EE'_J\to\FF$ lifts uniquely to $\EE_J$.
By cohomology and base change, this is also true over an arbitrary base $T$.
Thus, restriction gives an isomorphism
$$
\theta_J:
\Hom_{\PP^n}(\EE'_J,\EE_J)\times\QQ_J
\xrightarrow{\ \sim\ }\UU_J,
\qquad
(h,g)\longmapsto[g\circ(h,\mathrm{id})].
$$

Let $\mathbf G_J\subseteq\mathbf G$ be the standard chart where the coordinates indexed by $J$ map isomorphically onto the Grassmannian quotient and let $ P_J= \Flag(\nu;\PP^n)\times \PP_{\mathbf G_J}(\mathscr V_c|_{\mathbf G_J})$.
A quotient in $\mathbf G_J$ is uniquely represented by a matrix $(h,\mathrm{id})$, with identity columns indexed by $J$.
Since the only nonzero maps from $\EE'_J$ to $\EE_J$ are between summands of twist $-a_q$, this gives $\mathbf G_J\simeq\Hom_{\PP^n}(\EE'_J,\EE_J)$ and $P_J\simeq\Hom_{\PP^n}(\EE'_J,\EE_J)\times\QQ_J$.
Under this isomorphism, $\theta_J$ sends a flag and a polynomial vector to the quotient obtained by pulling back its rank-one kernel along $(h,\mathrm{id})$. 
Since this construction is independent of the choice of identity columns, the maps agree on overlaps and glue to a morphism
$$
\Theta: \Flag(\nu;\PP^n)\times\PP_{\mathbf G}(\mathscr V_c) \longrightarrow\QQ_{\bfa}(s,\mu).
$$
The isomorphisms on the charts glue to an isomorphism onto $\bigcup_J\UU_J$. 
This union is open and, since the source of $\Theta$ is projective, it is also closed.
Since $\QQ_{\bfa}(s,\mu)$ is connected, it follows that the union is the whole Quot scheme, as desired.
The case of $\mu = \varnothing$ is analogous.
\end{proof}

\begin{definition} \label{defn_residual}
Write $\lambda=(\mu,1^b)$, where $\mu$ has no parts equal to one.
The \textbf{residual quotient functor of type $(\bfa,s,\lambda)$} is
$$
\ResQuot_{\bfa}(s,\lambda): (\mathrm{Sch}/\kk)^{\mathrm{op}} \longrightarrow \mathrm{Sets},
$$
which assigns to a finite type $\kk$-scheme $T$ the set of equivalence classes of chains of coherent sheaves
$(\EE_{\bfa})_T \twoheadrightarrow \FF \twoheadrightarrow \FF_0$ such that
\begin{itemize}
\item$\FF_0$ is $T$-flat with Hilbert polynomial $P_{\bfa,s,\mu}$;
\item $\ker(\FF\to\FF_0)$ is $T$-flat with fibers of length $b$
\end{itemize}
Two chains are equivalent if they have the the same sequence of kernel subsheaves of $(\EE_{\bfa})_T$.
We call such a chain a \textbf{residual quotient flag}, and its
first quotient $(\EE_{\bfa})_T\twoheadrightarrow\FF$ the associated
\textbf{residual quotient}.

The \textbf{configuration functor} $\ResQuot_{\bfa}^{\circ}(s,\lambda)$ is the open subfunctor consisting
of flags for which $\ker(\FF\to\FF_0)$ is an invertible sheaf on a finite \'etale degree-$b$ subscheme contained in the locally free locus of $\FF_0$.
\end{definition}

\begin{remark} If $[\EE_{\bfa} \twoheadrightarrow\FF\twoheadrightarrow\FF_0]$ is a $\kk$-point of the configuration locus, then
$$
\ker(\FF\to\FF_0)\simeq\bigoplus_{i=1}^b\kk_{p_i}
$$
for distinct points $p_1,\ldots,p_b$ at which $\FF_0$ is locally free as an $\OO_{\PP^n}$-module. 
In the case of Hilbert schemes, $\FF_0=\OO_Y$, and the configuration locus therefore has, as part of its parameterization, $b$ distinct points disjoint from $Y$.
\end{remark}

\begin{example}
Let $\EE_{(0)}=\OO_{\PP^2}$ and $\lambda = (2,1)$. Let $L=V(x)\subseteq\PP^2$ and note that $\OO_L$ is locally free on $\PP^2 \setminus L$.

For $p \notin L$ the chain
$$
[\OO_{\PP^2}\twoheadrightarrow \OO_{L\sqcup\{p\}}\twoheadrightarrow\OO_L]
$$
is a residual quotient flag in $\ResQuot_{(0)}(0,\lambda)$ and also belongs to the configuration locus.

Let $p=[0:0:1] \in L$ and $Z=V(x^2,xy)$, and note that the chain
$$
[\OO_{\PP^2}\twoheadrightarrow\OO_Z\twoheadrightarrow\OO_L]
$$
is also a residual quotient flag. 
However, this time $\ker(\OO_Z \to \OO_L) = \kk_p$ and, since $p$ is not in the locally free locus, the flag does not lie in the configuration locus.
\end{example}

\begin{thm} \label{residual_quot} The residual quotient functor is represented by
$$\ResQuot_{\bfa}(s,\lambda) := \Quot^b_{\PP^n\times \QQ_{\bfa}(s,\mu)/\QQ_{\bfa}(s,\mu)}(\KK_\mu),$$
and forgetting the intermediate quotient defines a closed immersion
$\iota:\ResQuot_{\bfa}(s,\lambda) \hookrightarrow\QQ_{\bfa}(s,\lambda)$.
Moreover, the closure of the configuration locus is  the lexicographic component $\LL_{\bfa}(s,\lambda)$.
\end{thm}

\begin{proof}
Recall the universal sequence 
$$
0 \longrightarrow \KK_\mu \longrightarrow (\EE_\bfa)_{\QQ_{\bfa}(s,\mu)} \longrightarrow \FF_\mu \longrightarrow 0.
$$
By construction, the functor $\ResQuot_{\bfa}(s,\lambda)$ is the relative Quot functor $\Quot^b_{\PP^n\times \QQ_{\bfa}(s,\mu)/\QQ_{\bfa}(s,\mu)}(\KK_\mu)$; this is representable by a projective scheme \cite[Theorem~2.2.4]{HL}. 
The configuration locus is open because the \'etale locus and locally free locus are both open.
The morphism $\iota$ is defined on $T$-points by
$$
[(\KK_\mu)_T \twoheadrightarrow \mathscr{T}] \mapsto
\left[
(\EE_{\bfa})_T\twoheadrightarrow
\frac{(\EE_{\bfa})_T}
{\ker((\KK_\mu)_T\twoheadrightarrow \mathscr{T})}
\right].
$$
Every geometric fiber of $\FF_\mu$ has no nonzero zero-dimensional subsheaf.
Indeed, on a chart in \Cref{core_isom}, the fiber fits into an exact sequence
$$
0 \to \OO_Y(-d) \to (\FF_\mu)_z \to \EE_J/f(\OO_{\PP^n}(-d)) \to 0.
$$
The first term has no zero-dimensional sheaf since $Y$ is a residual flag while the last term has depth at least $n-1$; in particular, $\FF_\mu$ has no nonzero torsion subsheaf.
For a geometric point of the relative Quot scheme, we have an exact sequence
$$
0 \longrightarrow \mathscr{T} 
\longrightarrow (\EE_{\bfa})/\ker( (\KK_\mu)_z \twoheadrightarrow \mathscr{T}) 
\longrightarrow (\FF_\mu)_z \longrightarrow0.
$$
Since $\mathscr{T}$ has length $b$, it is precisely the maximal zero-dimensional subsheaf of the middle
term. Thus the resulting quotient uniquely determines the original residual quotient flag.
Hence $\iota$ is injective on geometric points. 
Its relative tangent space is $\Hom_{\PP^n}(\mathscr{T}, (\FF_\mu)_z)=0$, so $\iota$ is unramified and, since it is proper, it is a closed immersion \cite[Tag~04XV]{stacks_project}.

Let $U\subseteq \PP^n \times \QQ_{\bfa}(s,\mu)$ be the locally free locus of $\FF_\mu$. 
Over $U$, the kernel $\KK_\mu$ is locally free of rank $q$. 
A point of the configuration locus is therefore a choice of $b$ distinct points of $U$ (lying over the same point of $\QQ_{\bfa}(s,\mu)$), together
with a one-dimensional quotient of $\KK_\mu$ at each point, modulo permutation. 	
The ordered choices form an open subset of an iterated projective bundle over an open subset of $(\PP^n)^b\times \QQ_{\bfa}(s,\mu)$; this space is smooth and irreducible. 
The symmetric group acts freely by permuting the points and thus its finite \'etale quotient, which is the configuration locus, is also smooth and irreducible.

At a point of this locus lying over $z$, local freeness of $(\FF_\mu)_z$ at the support of $\mathscr{T}$ and $n \geq 2$ give $\Hom_{\PP^n}(\mathscr{T},\FF_\mu)=\Ext^1_{\PP^n}(\mathscr{T}, (\FF_\mu)_z)=0$.
By \Cref{Deformation_theory_flag}, $\iota$ is \'etale there. Its restriction to the configuration locus is therefore an open immersion.

By \Cref{hilbert_lex}, $L(\lambda)$ is a specialization of a residual flag of type $\mu$ together with $b$ distinct points outside it.
Putting these ideals in the $q$th summand of $E_{\bfa}$, adjoining the first $q-1$ full summands and the last $s$ zero summands, gives a specialization from $\iota(\ResQuot_{\bfa}^{\circ}(s,\lambda))$ to the lexicographic point of $\QQ_{\bfa}(s,\lambda)$.
Since the configuration locus is irreducible and its image is open, its closure is an irreducible component. 
It contains the smooth lexicographic point by \Cref{lex_smooth}, so it is $\LL_{\bfa}(s,\lambda)$.
\end{proof}

\begin{cor} \label{lex_general_member}
The general member of $\LL_{\bfa}(s,\lambda)$ is obtained from a general quotient in $\QQ_{\bfa}(s,\mu)$
by taking length-one quotients of its kernel at $b$ distinct general points of $\PP^n$. 
Moreover, $\LL_{\bfa}(s,\lambda)$ is generically smooth of dimension $\dim\QQ_{\bfa}(s,\mu)+b(n+q-1)$.
\end{cor}
\begin{proof}
We only need to check the dimension. 
For each of the $b$ points $p$, we first choose its position in $\PP^n$ and then a quotient
$(\KK_\mu)_p \otimes \kk(p) \twoheadrightarrow \kk(p)$. 
Since $\KK_\mu$ is locally free of rank $q$ at these points, quotients are parametrized by $\PP^{q-1}$. 
Thus $\dim\LL_{\bfa}(s,\lambda) =\dim\QQ_{\bfa}(s,\mu)+b(n+q-1)$.
\end{proof}

\subsection{Reeves' double saturation}
We now generalize Reeves' criterion for a Borel-fixed point to lie on the lexicographic component \cite{Reeves_radius}.

\begin{definition} \label{double_sat}
The \textbf{double saturation} of a monomial submodule $K\subseteq E_{\bfa}$ is
$$
\delta(K)=(K:x_n^\infty):x_{n-1}^\infty = K:(x_{n-1}x_n)^\infty.
$$
\end{definition}

\begin{remark}
The double saturation is obtained by setting $x_{n-1}=x_n=1$ in the monomial generators. 
Since our submodules are saturated, if $K$ is Borel-fixed, then $\delta(K)=K:x_{n-1}^\infty$.
In particular, if $\lambda=(\mu,1^b)$ with $\mu$ having no parts equal to one, then $\delta(L_{\bfa}(s,\lambda))=L_{\bfa}(s,\mu)$.
\end{remark}

We need a simple extension of \cite[Lemma~2.1]{SSS26}.

\begin{lemma}\label{splitting}
Let $\lambda,\nu_1,\ldots,\nu_m$ be partitions and $b_1,\ldots,b_m$ nonnegative integers such that
$$
p_\lambda(t)=\sum_{i=1}^m p_{\nu_i}(t+b_i).
$$
If $\lambda$ has no parts equal to one, then at most one of the $\nu_i$ is nonempty.
\end{lemma}

\begin{proof}
Suppose $\nu_1\ne\varnothing$, and let $Z_i$ be the scheme defined by $L(\nu_i)$. 
The sheaves $\OO_{Z_1}(b_1)$ and $\bigoplus_{i=2}^m\OO_{Z_i}(b_i)$ are globally generated and, in particular, are quotients of trivial bundles. 
Since both have rank zero, \Cref{Macaulay}, applied with all twists equal to zero, expresses their Hilbert polynomials as $p_\nu$ and $p_\mu$ for some partitions $\nu$ and $\mu$. 
The first sheaf is nonzero, so $\nu\ne\varnothing$. 
Since $p_\lambda=p_\nu+p_\mu$, \cite[Lemma~2.1]{SSS26} implies that $\mu$ is empty, and thus $\nu_i=\varnothing$ for every $i\geq2$.
\end{proof}

\begin{lemma}\label{residual_unique}
If $\lambda$ has no parts equal to one, then $L_{\bfa}(s,\lambda)$ is the unique Borel-fixed point of $\QQ_{\bfa}(s,\lambda)$.
\end{lemma}

\begin{proof}
By \Cref{Pardue_Borel}, a Borel-fixed kernel is of the form $K=\bigoplus_{i=1}^r I_i(-a_i)\bfe_i$, with the $I_i$ saturated Borel-fixed ideals satisfying $I_1\supseteq\cdots\supseteq I_r$. 
Since $K$ has rank $q$, we have $I_i=0$ for $i>q$. 
For $i \leq q$, each $[S/I_i]$ is a point on a Hilbert scheme $\Hilb^{p_{\nu_i}}(\PP^n)$ for some partition $\nu_i$.
Comparing Hilbert polynomials we have
$$
p_\lambda(t)=\sum_{i=1}^q p_{\nu_i}(t+a_q-a_i).
$$
But, \Cref{splitting} forces at most one of the $\nu_i$ to be non-zero, which, by nesting,  force $I_i=S$ for $i<q$ and $P_{S/I_q}=p_\lambda$. 
By \cite[Theorem~1.1]{Staal_smooth_hilbert_schemes}, the only saturated Borel-fixed ideal with this polynomial is $L(\lambda)$.
Thus $K=L_{\bfa}(s,\lambda)$.
\end{proof}

\begin{thm}\label{double_saturation}
Assume $n\geq2$, and let $[E_{\bfa}/K]$ be a Borel-fixed point of $\QQ_{\bfa}(s,\lambda)$. 
Then $[E_{\bfa}/K]$ lies on the lexicographic component if and only if $K$ and $L_{\bfa}(s,\lambda)$ have the same double saturation.
\end{thm}

\begin{proof}
Write $\lambda=(\mu,1^b)$, with $\mu$ having no parts equal to one, and put $L_0=L_{\bfa}(s,\mu)$.
The case $b=0$ follows from \Cref{residual_unique}, so assume $b>0$.

Suppose first that $[E_{\bfa}/K]\in\LL_{\bfa}(s,\lambda)$.
By \Cref{lex_general_member}, this quotient has a unique residual quotient flag 
$[E_{\bfa} \twoheadrightarrow E_{\bfa}/K \twoheadrightarrow F_0]$ with $[E_{\bfa} \twoheadrightarrow F_0] \in \QQ_{\bfa}(s,\mu)$.
Since $\ker(E_{\bfa} \twoheadrightarrow F_0)$ is also Borel-fixed,  \Cref{residual_unique} implies that this kernel must be $L_{\bfa}(s,\mu)$. 
Thus, $K\subseteq L_{\bfa}(s,\mu)$ and the sheaf $\widetilde{L_{\bfa}(s,\mu)/K}$ has length $b$.
Its support is the Borel-fixed point $p=[0:\cdots:0:1]$, so a power of $x_{n-1}$ annihilates this sheaf.
Since $K$ is saturated, it follows that
$$
L_{\bfa}(s,\mu) \subseteq K:x_{n-1}^\infty \subseteq L_{\bfa}(s,\mu):x_{n-1}^\infty = L_{\bfa}(s,\mu).
$$
Thus $\delta(K) = L_{\bfa}(s,\mu) = \delta(L_{\bfa}(s,\lambda))$, as required.

Conversely, suppose $\delta(K)=L_{\bfa}(s,\mu)$, and let $K=\bigoplus_i I_i(-a_i)\bfe_i$.
Since $K \subseteq L_{\bfa}(s,\mu)$, comparing Hilbert polynomials we see that $\widetilde{L_{\bfa}(s,\mu)/K}$ has length $b$.
It follows that for $i<q$, the ideal $I_i$ defines a scheme of length $b_i$, the sheaf $\widetilde{L(\mu)/I_q}$ has length $b_q$, and $I_i=0$ for $i>q$. 
We also have  $b_1 + \cdots + b_q=b$.

By \Cref{lemma_distraction}, we can choose compatible distraction parameters for $K\subseteq L_{\bfa}(s,\mu)$. 
In particular, if we let $J_i(t)=D_t(I_i)$ and $J_\mu(t)=D_t(L(\mu))$ then, by compatibility, $J_q(t)\subseteq J_\mu(t)$ for all $t$, and these ideals agree with $I_q$ and $L(\mu)$, respectively, at $t=0$. 
The family $\gamma_K$ of \Cref{lemma_distraction} is given by $\gamma_K(t)=[E_{\bfa}/K_t]$ with $K_t=\bigoplus_{i=1}^q J_i(t)(-a_i)\bfe_i$.
Let $L_t = \bigoplus_{i<q}S(-a_i)\bfe_i \oplus J_\mu(t)(-a_q)\bfe_q$ and we obtain a family of
residual quotient flags
$$
E_{\bfa}\twoheadrightarrow E_{\bfa}/K_t \twoheadrightarrow E_{\bfa}/L_t.
$$
For $t\ne0$, the distraction parameters can be chosen so that all $J_i(t)$ and $J_\mu(t)$ are radical, the supports $V(J_i(t))$ for $i<q$ are pairwise disjoint and avoid $V(J_q(t))$, and $V(J_q(t))$ is the disjoint union of $V(J_\mu(t))$ and $b_q$ reduced points.
The kernel of the second quotient is therefore a sum of $b$ skyscraper sheaves at distinct points.
Thus, the residual quotient flag belongs to the configuration locus for every $t\ne0$. By \Cref{lex_general_member}, the special fiber $[E_\bfa/K]$ lies in the closure of the configuration locus, namely $\LL_{\bfa}(s,\lambda)$, as required.
\end{proof}

\section{Smooth Quot schemes}
We classify the smooth schemes $\QQ_{\bfa}(s,\lambda)$ and use the residual quotient construction to describe them. 
When all the $a_i$ are equal, the classification specializes to \cite[Theorem~1.1]{SSS26}.
We end the section by proving that the Quot scheme of the tangent bundle to $\PP^2$ is connected.

\begin{lemma}\label{n1Borel}
Suppose $n=1$ and $\lambda=(1^b)$. 
All the Borel-fixed points of $\QQ_{\bfa}(s,\lambda)$ lie on a common irreducible component and are of the form $[E_{\bfa}/K_{\boldsymbol\beta}]$, where
$$
K_{\boldsymbol\beta}
=
\bigoplus_{i=1}^q x_0^{\beta_i}S(-a_i)\bfe_i \oplus \bigoplus_{j>q}0\bfe_j,
$$
for some $ 0\leq\beta_1\leq\cdots\leq\beta_q$ with $\sum_{i=1}^q\beta_i=b$.
Their tangent spaces have dimensions
$$
\dim_\kk T_{[E_{\bfa}/K_{\boldsymbol\beta}]}\QQ_{\bfa}(s,\lambda)
= qb+\sum_{i=1}^q\sum_{j>q} \max\{a_i+\beta_i-a_j+1,0\}.
$$
\end{lemma}

\begin{proof}
The description of the Borel-fixed points follows from \Cref{Pardue_Borel} and the fact that they all lie on one component was shown in \Cref{radius_bound}.
Let $B_i = V(x_0^{\beta_i})$ and note that we have 
$$
\widetilde{K_{\boldsymbol\beta}} \simeq\bigoplus_{i=1}^q\OO_{\PP^1}(-a_i-\beta_i),
\,\, \text{ and } \,\,
\widetilde{E_{\bfa}/K_{\boldsymbol\beta}} \simeq \bigoplus_{j=1}^q\OO_{B_j}(-a_j) \oplus\bigoplus_{j>q}\OO_{\PP^1}(-a_j).
$$
Thus the tangent space to $[E_{\bfa}/K]$ is
$$
T_{[E_{\bfa}/K]}\QQ_{\bfa}(s,\lambda)
= \Hom_{\PP^1}(\widetilde{K_{\boldsymbol\beta}}, \widetilde{E_{\bfa}/K_{\boldsymbol\beta}} )
=
\bigoplus_{i=1}^q
H^0\bigl(\PP^1, \widetilde{E_{\bfa}/K_{\boldsymbol\beta}} (a_i+\beta_i)\bigr).
$$
For each $i$, the torsion summands contribute $\sum_{j=1}^q\beta_j=b$ dimensions, since twisting a zero-dimensional sheaf does not change the dimension of its sections. 
The remaining summands contribute
$$
\sum_{j>q}h^0\bigl(\OO_{\PP^1}(a_i+\beta_i-a_j)\bigr) = \sum_{j>q}\max\{a_i+\beta_i-a_j+1,0\}.
$$
Summing over all $i$ gives the desired formula.
\end{proof} 

To prove that certain Quot schemes are singular when $n\geq 2$, we need to find a Borel-fixed point on the lexicographic component and prove that the tangent space dimension is large. 
The following lemma establishes some of the numerics of this other Borel-fixed point.

\begin{lemma} \label{tan_excess}
Let $n\geq2$, $q\geq2$, and let $\lambda=(\nu,1)$ with $\nu$ non-empty. 
Let $p=[0:\cdots:0:1]$, $\delta=a_q-a_{q-1}$, and $Z_\nu=V(L(\nu))$.
Consider the Borel-fixed kernel
$$
K'=\bigoplus_{i<q-1}S(-a_i)\bfe_i \oplus I_p(-a_{q-1})\bfe_{q-1} \oplus L(\nu)(-a_q)\bfe_q \oplus\bigoplus_{j>q}0\bfe_j.
$$
Then we have 
\begin{align}\label{eq_tan_ex}
&\dim_\kk T_{[E_{\bfa}/K']}\QQ_{\bfa}(s,\lambda)
-\dim\LL_{\bfa}(s,\lambda) \notag\\
&=n-\nu_{\ell(\nu)}
+\dim_\kk\Hom_{\PP^n}(\II_p,\OO_{Z_\nu}(-\delta)) - h^0(\OO_{Z_\nu}(-\delta)).
\end{align}
\end{lemma}
\begin{proof}
For any partition $\rho$, let $Z_\rho=V(L(\rho))$ and let $d_\rho$ denote the dimension of its lexicographic component on $\Hilb^{p_\rho}(\PP^n)$.
Let $\KK'$ and $\FF'$ denote the sheafified kernel and quotient at $K'$, and let
$\KK_{\lex}$ and $\FF_{\lex}$
denote those at the lexicographic point $[E_{\bfa}/L_{\bfa}(s,\lambda)]$. 
Explicitly,
$$
\mathscr K'= \bigoplus_{i<q-1}\OO_{\PP^n}(-a_i) \oplus\II_p(-a_{q-1}) \oplus\II_{Z_\nu}(-a_q), \, \text{and} \,
\mathscr F' = \OO_p(-a_{q-1}) \oplus\OO_{Z_\nu}(-a_q) \oplus\bigoplus_{j>q}\OO_{\PP^n}(-a_j).
$$
Similarly,
$$
\KK_{\lex}  = \bigoplus_{i<q}\OO_{\PP^n}(-a_i) \oplus\II_{Z_\lambda}(-a_q) \,\, \text{and} \, \, 
\mathscr F_{\lex} = \OO_{Z_\lambda}(-a_q) \oplus\bigoplus_{j>q}\OO_{\PP^n}(-a_j).
$$
Since $\II_p^\vee\simeq\OO_{\PP^n}$ and $\II_{Z_\nu}^\vee\simeq\II_{Z_\lambda}^\vee$, we have $(\KK')^\vee\simeq\KK_{\lex}^\vee$. 
Hence, for every $j>q$,
$$
\Hom_{\PP^n}(\KK',\OO_{\PP^n}(-a_j)) \cong H^0(\PP^n,(\mathscr K')^\vee(-a_j))
\cong H^0(\PP^n,\mathscr K_{\lex}^\vee(-a_j)) \cong \Hom_{\PP^n}(\KK_{\lex},\OO_{\PP^n}(-a_j)).
$$
The terms with target $\OO_{\PP^n}(-a_j)$ therefore cancel in the difference below. 
Expanding the remaining Hom spaces along the direct sums, we obtain
\begin{align*}
&\dim_\kk\Hom_{\PP^n}(\KK',\FF') -\dim_\kk\Hom_{\PP^n}(\KK_{\lex}, \FF_{\lex})\\
&=\sum_{i<q-1}\Big(
\dim_\kk\Hom_{\PP^n}(\OO_{\PP^n}(-a_i),\OO_p(-a_{q-1}))
+\dim_\kk\Hom_{\PP^n}(\OO_{\PP^n}(-a_i),\OO_{Z_\nu}(-a_q))\\
&\hspace{42mm}
-\dim_\kk\Hom_{\PP^n}(\OO_{\PP^n}(-a_i),\OO_{Z_\lambda}(-a_q))
\Big)\\
&\quad + \dim_\kk\Hom_{\PP^n}(\II_p,\OO_p)
+\dim_\kk\Hom_{\PP^n}(\II_p,\OO_{Z_\nu}(-\delta))\\
&\quad + \dim_\kk\Hom_{\PP^n}(\II_{Z_\nu},\OO_p(\delta))
+\dim_\kk\Hom_{\PP^n}(\II_{Z_\nu},\OO_{Z_\nu})
\\
&\quad-\dim_\kk\Hom_{\PP^n}(\OO_{\PP^n}(-a_{q-1}),\OO_{Z_\lambda}(-a_q))
-\dim_\kk\Hom_{\PP^n}(\II_{Z_\lambda},\OO_{Z_\lambda})
\end{align*}
Using the exact sequence
$$
0\longrightarrow\OO_p\longrightarrow \OO_{Z_\lambda}\longrightarrow\OO_{Z_\nu} \longrightarrow0.
$$
and twisting it by $a_i-a_q$, we get $ h^0(\OO_{Z_\lambda}(a_i-a_q)) = h^0(\OO_{Z_\nu}(a_i-a_q))+1$. 
It follows that the sum in the parantheses is zero. 
We also have 
$$
\dim_\kk\Hom_{\PP^n}(\II_{Z_\nu},\OO_p(\delta)) = \dim_\kk \Hom_\kk(\II_{Z_\nu}/\mfm_p\II_{Z_\nu}, \kk_p)
$$
which is precisely the number of minimal generators of $\II_{Z_\nu}$ (\Cref{mingen}), which is $n- \nu_{\ell(\nu)}+1$.
We have $\dim_\kk \Hom_{\PP^n}(\II_p,\OO_p) = n$, since it is the tangent space to $p \in \PP^n$ on its Hilbert scheme.
By \cite[Theorem 4.1]{RS97} and the smoothness of the lexicographic point we obtain
$$
\dim_\kk\Hom_{\PP^n}(\II_{Z_\lambda},\OO_{Z_\lambda}) = d_\lambda = d_\nu+n 
= \dim_\kk\Hom_{\PP^n}(\II_{Z_\nu},\OO_{Z_\nu}) + n
$$
Since $[\EE_{\bfa} \to \FF_{\lex}]$ is a smooth point on the lexicographic component, combining all these dimensions gives the desired result.
\end{proof}

Finally, we need two more simple lemmas.

\begin{lemma} \label{hom_max}
Let $(R,\mfm,\kk)$ be an $n$-dimensional regular local ring with a regular sequence $y_1,\dots,y_n$. 
Let $M = R/(y_1,\dots,y_{n-1},y_n^\ell)$ for some $\ell \geq 1$. 
Then we have  $\dim_\kk \Hom_R(\mfm,M)=\ell+n-1$.
\end{lemma}
\begin{proof}
The syzygies of $\mfm$ allow us to identify the space as  $\Hom_R(\mfm,M) =\{(m_1,\ldots,m_n)\in M^n: y_im_j=y_jm_i\}$.
Since $y_1,\ldots,y_{n-1}$ annihilate $M$, the only condition is $y_nm_i = 0$ for $i \leq n-1$. 
In particular, the Hom space is just 
$$
(y_{n}^{\ell-1}/y_n^{\ell})^{\oplus(n-1)} \oplus M.
$$
The dimension follows immediately. 
\end{proof}

\begin{lemma}\label{hom_hypersurface}
Let $p=[0:0:1]\in\PP^2$ and $Z = V(x_0^k)$ for any $k\geq 1$.
For every $\delta \geq 0$ we have
$$
\dim_\kk \Hom_{\PP^2}(\II_p,\OO_{Z}(-\delta)) - h^0(\OO_{Z}(-\delta))
=
\begin{cases}
1, & k + \delta \leq 2,\\
0, & k + \delta \geq 3.
\end{cases}
$$
\end{lemma}

\begin{proof}
Applying $\HH om_{\PP^2}(-,\OO_Z)$ to the Koszul resolution
$$
0\longrightarrow\OO_{\PP^2}(-2)
\xrightarrow{(-x_1,x_0)}
\OO_{\PP^2}(-1)^{\oplus2}
\longrightarrow\II_p\longrightarrow0
$$
identifies $\mathscr H:=\mathcal Hom_{\PP^2}(\II_p,\OO_Z)$ with $\ker\bigl(\OO_Z(1)^{\oplus2} \xrightarrow{(x_1,-x_0)}\OO_Z(2)\bigr)$.
Since $Z=V(x_0^k)$, we can resolve this as follows
$$
0\longrightarrow\OO_{\PP^2}(1-k)^{\oplus2}
\xrightarrow{
\begin{pmatrix}
x_0^{k-1}&0\\
-x_1&x_0
\end{pmatrix}}
\OO_{\PP^2}\oplus\OO_{\PP^2}(2-k)
\longrightarrow\mathscr H\longrightarrow0.
$$
Taking cohomology and using the fact that $k\geq 1$ and $\delta \geq 0$, we obtain
\begin{align*}
h^0(\mathscr H(-\delta))  & = h^0(\OO_{\PP^2}(-\delta)) +h^0(\OO_{\PP^2}(2-k-\delta)) -2h^0(\OO_{\PP^2}(1-k-\delta)) \\
& = h^0(\OO_{\PP^2}(-\delta)) +h^0(\OO_{\PP^2}(2-k-\delta)).
\end{align*}
On the other hand, the ideal sheaf sequence for $Z$ gives  $h^0(\OO_Z(-\delta))=h^0(\OO_{\PP^2}(-\delta))$.
Subtracting these two, gives the desired result.
\end{proof}

We are now ready for the main theorem.

\begin{thm}\label{smooth_quot_classification}
The Quot scheme $\QQ_{\bfa}(s,\lambda)$ is smooth if and only if one of the following holds.
\begin{enumerate}
\item $n=1$, $\lambda=(1^b)$, and either $b=0$, $q=1$, or there is no
$j>q$ such that $a_1+2\leq a_j\leq a_q+b$.
\item $n=2$ and $q=1$.
\item $n\geq3$, $q=1$, and, writing $\lambda=(n^c,\nu)$ with $\nu$ having no parts equal to $n$, one of the following holds:
\begin{enumerate}
\item $\nu$ has no parts equal to one, including $\nu=\varnothing$;
\item $\nu=(1^d)$ with $1\leq d\leq3$;
\item $\nu=(m^d,1)$ with $2\leq m<n$ and $d\geq1$, excluding
$(m,d)=(2,2)$ and $(2,3)$.
\end{enumerate}
\item $n\geq2$, $q\geq2$, and $\lambda = (1)$ or has no parts equal to one.
\item $n\geq2$, $q\geq2$, and $\lambda=(n^c,1)$ with $c\geq1$, where either $n\geq3$, or $n=2$ and $c+a_q-a_{q-1}\geq3$.
\end{enumerate}
In each case, the Quot scheme is irreducible.
\end{thm}

\begin{proof}
By \Cref{Pardue_conn}, a smooth Quot scheme is also irreducible. 
Since the singular locus is closed and invariant under $B_{\bfa}$, it suffices to determine when every Borel-fixed point on the Quot scheme is smooth. 

Assume $n=1$. If $b=0$ or $q=1$, the lexicographic point is the unique fixed point and is smooth by \Cref{lex_smooth_all}. 
If $b>0$ and $q\geq2$, \Cref{lex_smooth_all} also shows that the conditions in (i) are necessary. 
To see that they are sufficient note that by \Cref{n1Borel} for every Borel-fixed point $[E_{\bfa}/K_{\boldsymbol\beta}]$ we have
$$
\dim_\kk T_{[E_{\bfa}/K_{\boldsymbol\beta}]}\QQ_{\bfa}(s,\lambda) = qb+\sum_{i=1}^q\sum_{j>q} \max\{a_i+\beta_i-a_j+1,0\}.
$$
Under the hypothesis in (i), each $j>q$ satisfies either $a_j\leq a_1+1$ or $a_j>a_q+b$. In the first case,
$$
\sum_{i=1}^q\max\{a_i+\beta_i-a_j+1,0\}
=\sum_{i=1}^q(a_i+\beta_i-a_j+1)
=\sum_{i=1}^q(a_i-a_j+1)+b,
$$
while in the second case, every summand is zero.
Thus the tangent dimension is independent of $[E_{\bfa}/K_{\boldsymbol\beta}]$ and equals the dimension at the smooth lexicographic point.
By \Cref{n1Borel}, all the Borel-fixed points lie on a common irreducible component which, in this case, must be the lexicographic component. 
Thus, all the Borel-fixed points are smooth as required.

Now suppose $n\geq2$ and $q=1$. By \Cref{P1_trivial}, smoothness is equivalent to that of $\Hilb^{p_\nu}(\PP^n)$, where $\lambda=(n^c,\nu)$ and $\nu$ has no parts equal to $n$.
For $n=2$, this is a Hilbert scheme of points on a smooth surface.
For $n\geq3$, \cite[Theorem~A]{SS23} gives precisely (iii).

Assume $n\geq2$ and $q\geq2$.
The case when $\lambda$ has no parts equal to one follows from \Cref{core_isom}.
If $\lambda=(1)$, the nesting forces the length-one quotient to occur in the $q$th summand; thus, the lexicographic point is the unique Borel-fixed point, and the result follows from \Cref{lex_smooth}.

So we may assume that $\lambda=(\nu,1)$ with $\nu \ne \varnothing$.
Let $p=[0:\cdots:0:1]$ and consider the Borel-fixed kernel
$$
K'=\bigoplus_{i<q-1}S(-a_i)\bfe_i \oplus I_p(-a_{q-1})\bfe_{q-1} \oplus L(\nu)(-a_q)\bfe_q \oplus\bigoplus_{j>q}0\bfe_j.
$$
The Hilbert polynomial of its quotient is $P_{\bfa,s,\lambda}$ since $p_\nu+1=p_\lambda$ and its double saturation agrees with that of $L_{\bfa}(s,\lambda)$. 
It follows from  \Cref{double_saturation} that $[E_{\bfa}/K']$ also lies on $\LL_{\bfa}(s,\lambda)$.

Assume further that $\lambda$ has at least two parts equal to one, i.e., $\nu=(\mu,1^\ell)$ with $\ell\geq1$ and $\mu$ having no parts equal to one. 
Consider the exact sequence
\begin{align} \label{Torsion_exact}
0 \longrightarrow \mathscr{T} \longrightarrow \OO_{Z_\nu} \longrightarrow \OO_{Z_\mu} \longrightarrow 0.
\end{align}
On the open chart $D_{+}(x_n)$, the module defining $\mathscr{T}|_{D_{+}(x_n)}$ is precisely $\kk[x_{n-1}]/(x_{n-1}^\ell)$.
By \Cref{hom_max}, we obtain $\dim_\kk\Hom_{\PP^n}(\II_p,\OO_{Z_{\nu}}(-\delta)) \geq \dim_\kk\Hom_{\PP^n}(\II_p,\mathscr{T}(-\delta)) = \ell+n-1$.
Taking cohomology of \Cref{Torsion_exact} and using \Cref{negative_cohomology} we obtain
$h^0(\OO_{Z_\nu}(-\delta))\leq\ell+1$.
Consequently, \Cref{eq_tan_ex}
gives
$$
\dim_\kk T_{[E_{\bfa}/K']}\QQ_{\bfa}(s,\lambda) - \dim\LL_{\bfa}(s,\lambda) \geq 2n-3>0.
$$
Since $[E_{\bfa}/K']$ lies on the lexicographic component, it must be a singular point.

Thus, we may assume that $\lambda$ has exactly one part equal to one. 
Then $Z_\nu$ has no zero-dimensional subsheaf, which implies that $\HH om_{\PP^n}(\OO_p,\OO_{Z_{\nu}})= 0$.
Applying $\HH om_{\PP^n}(-,\OO_{Z_{\nu}})$ to the ideal sheaf sequence of the point $p\in \PP^n$ we obtain an inclusion $\OO_{Z_\nu}\hookrightarrow \mathcal{H}om_{\PP^n}(\II_p,\OO_{Z_\nu})$.
Hence $\dim_\kk\Hom_{\PP^n}(\II_p,\OO_{Z_{\nu}}(-\delta)) \geq h^0(\OO_{Z_\nu}(-\delta))$. 
It follows from \Cref{eq_tan_ex} that $[E_\bfa/K']$ is singular as long as $\nu \ne (n^c)$.

Thus, we may assume $\nu=(n^c)$ with $c\geq1$. 
Since $L(\nu)=(x_0^c)$ and $\II_p^{\vee\vee}=\OO_{\PP^n}$, we have
$$
(\widetilde K')^{\vee\vee} = \bigoplus_{i<q-1}\OO_{\PP^n}(-a_i) \oplus \OO_{\PP^n}(-a_{q-1}) \oplus \OO_{\PP^n}(-a_q-c).
$$
This bundle is rigid, since $H^1(\PP^n,\OO_{\PP^n}(v))=0$ for every $v$.
Moreover, $(\widetilde K')^{\vee\vee}/\widetilde K' \simeq \OO_p(-a_{q-1})$.
Thus, for $n\geq3$, we can use \Cref{double_dual_cor} to reduce smoothness at $[E_{\bfa}/K']$ to smoothness of a length-one Quot scheme of a vector bundle, which is smooth; this follows, for example, by the canonical identification $\Quot^{1}_X(\EE) \simeq \PP_X(\EE^{\vee})$.
If $n=2$, using \Cref{hom_hypersurface} together with \Cref{tan_excess}, we see that $[E_\bfa/K']$ is singular iff $c+ (a_q - a_{q-1}) \leq 2$; this are precisely the cases excluded in (v). 
To finish, it will suffice to show that $[E_\bfa/K']$ and the lexicographic point are the only Borel-fixed points on the Quot scheme when $\lambda = (n^c,1)$ and satisfies the conditions in (v).

To this end, consider a Borel-fixed kernel $K = \bigoplus_i I_i(-a_i)\bfe_i$, with $I_1\supseteq\cdots\supseteq I_q\ne0$ and $I_i=0$ for $i>q$.
Comparing Hilbert polynomials we have  $p_\lambda(t)=\sum_{i=1}^q P_{S/I_i}(t+a_q-a_i)$.
Each nonzero term is the polynomial of a globally generated torsion sheaf, so \Cref{Macaulay} expresses it as $p_\rho$ for a nonempty partition.
By \cite[Lemma~2.1]{SSS26}, the only nontrivial decompositions are $p_{(n^c,1)}=p_{(n^c)}+1$ and, when  $n=c=2$, we also have $p_{(2,2,1)}=p_{(2)}+p_{(2)}$.
Neither decomposition admits a further nontrivial splitting.
If the splitting is just $p_{(n^c,1)}$ with one nonzero term, then the nesting forces the non-zero quotient to occur in position $q$ which in turn implies that the point is just the lexicographic point. 
If the splitting is $p_{(n^c)}+1$, the non-zero quotients occur in positions $q-1$ and $q$, which forces the kernel to be $K'$. 
Finally, if we have the splitting $p_{(2)} + p_{(2)}$ then $S/I_{q-1}$ would have Hilbert polynomial $t+1-\delta$, where $\delta=a_q-a_{q-1} \geq 0$. 
This implies that $a_q - a_{q-1} = 0$, which in turn implies $c+a_q - a_{q-1} =2$, contradicting the assumption of (v).
\end{proof}

\begin{cor}\label{smooth_quot_functor}
Suppose $n\geq2$ and write $\lambda=(\mu,1^b)$, where $\mu$ has no parts equal to one. If $\QQ_{\bfa}(s,\lambda)$ is smooth, then
$$
\QQ_{\bfa}(s,\lambda) \simeq \ResQuot_{\bfa}(s,\lambda) =\Quot^b_{\PP^n\times\QQ_{\bfa}(s,\mu)/\QQ_{\bfa}(s,\mu)}(\KK_\mu).
$$
\end{cor}
\begin{proof}
When the Quot scheme is smooth and irreducible, the closed immersion of \Cref{lex_general_member} becomes an isomorphism.
\end{proof}

We now describe the smooth Quot schemes in \Cref{smooth_quot_classification} for $n\geq 2$, in the same order. Recall the following notation.

\begin{notation}
The number of occurrences of $a_q$ on either side of the index $q$ are denoted by $q_0=\#\{i\leq q:a_i=a_q\}$ and $s_0=\#\{j>q:a_j=a_q\}$. 

Consider the Grassmannian $\mathbf{G} :=\Gr(s_0+1,q_0+s_0)$ and, for each $c \geq 0$, a locally free sheaf
$$
\mathscr{V}_c= (S_c\otimes_\kk \mathscr{Q}) \oplus \bigoplus_{a_j>a_q} S_{a_q+c-a_j}\otimes_\kk\OO_{\mathbf{G}}
$$
on it, where $\mathscr{Q}$ denotes the universal quotient.
\end{notation}

\begin{thm} \label{smooth_quot_geometry}
For $n\geq 2$, the smooth Quot schemes, $\QQ_{\bfa}(s,\lambda)$ have the following geometric descriptions.
\begin{enumerate}
\item
Suppose $n=2$ and $q=1$. If $\lambda=(2^c,1^b)$, then we have
$$
\QQ_{\bfa}(r-1,\lambda) \simeq \Hilb^b(\PP^2)\times\PP(\mathscr V_c).
$$
\item
Suppose $n\geq3$ and $q=1$.
If $\lambda=(n^c,\nu)$ with $\nu$ having no parts equal to $n$, then we have 
\begin{enumerate}
\item[(a)]
If $\nu$ has no parts equal to one, then
$$
\QQ_{\bfa}(r-1,\lambda) \simeq \Flag(\nu;\PP^n)\times\PP(\mathscr V_c).
$$
\item[(b)]
If $\nu=(1^d)$ with $1\leq d\leq3$, then
$$
\QQ_{\bfa}(r-1,\lambda) \simeq \Hilb^d(\PP^n)\times\PP(\mathscr V_c).
$$
\item[(c)]
If $\nu=(m^d,1)$ as in \Cref{smooth_quot_classification} (iii)(c), let $H=\Flag((m^d);\PP^n)$ and let $\mathcal Y\subseteq\PP^n\times H$ be the universal subscheme. 
Then
$$ \QQ_{\bfa}(r-1,\lambda) \simeq \Bl_{\mathcal Y}(\PP^n\times H) \times\PP(\mathscr V_c).$$
\end{enumerate}
\item
Suppose $n\geq2$ and $q\geq2$.
\begin{enumerate}
\item[(a)]If $\lambda=\varnothing$, then
$$
\QQ_{\bfa}(s,\varnothing) \simeq \Gr(s_0,q_0+s_0).
$$
\item[(b)] If $\lambda=(n^c,\nu)$ with $\nu$ has no parts equal to $n$ or $1$, then
$$
\QQ_{\bfa}(s,\lambda) \simeq \Flag(\nu;\PP^n) \times \PP_{\mathbf{G}}(\mathscr{V}_c).
$$
\item[(c)] If $\lambda=(1)$, let $B=\Gr(s_0,q_0+s_0)$.
Under the identification $B\simeq\QQ_{\bfa}(s,\varnothing)$, we have
$$
\QQ_{\bfa}(s,(1)) \simeq \PP_{\PP^n\times B}(\KK_\varnothing^\vee).
$$
In particular, this is a $\PP^{q-1}$-bundle over $\PP^n\times B$.
\end{enumerate}
\item
Suppose $\lambda=(n^c,1)$ as in \Cref{smooth_quot_classification} (v).
Let $B=\PP_{\mathbf G}(\mathscr V_c) \simeq\QQ_{\bfa}(s,(n^c))$. Then
$$
\QQ_{\bfa}(s,(n^c,1)) \simeq \PP_{\PP^n\times B}(\KK_{(n^c)}^\vee).
$$
In particular, this is a $\PP^{q-1}$-bundle over $\PP^n\times\PP_{\mathbf G}(\mathscr V_c)$.
\end{enumerate}
\end{thm}
\begin{proof}
Items (i), (ii)(a), and (ii)(b) follow from \Cref{P1_trivial}.

For item (ii)(c), let $\lambda=(\mu,1)$ with $\mu=(n^c,m^d)$.
\Cref{P1_trivial} once again gives $\QQ_{\bfa}(r-1,\mu)\simeq H\times\PP(\mathscr V_c)$.
Its universal kernel is the pullback of $\II_{\mathcal Y}$ tensored with an invertible sheaf.
Since tensoring by an invertible sheaf does not change the length-one Quot scheme, the residual quotient construction gives
$$
\QQ_{\bfa}(r-1,\lambda) \simeq \Proj_{\PP^n\times H}(\Sym\II_{\YY}) \times\PP(\mathscr V_c).
$$
The universal subscheme $\YY$ is a relative hypersurface in a relative linear subspace, hence is regularly embedded in $\PP^n\times H$.
By \cite[Proposition~2.2]{Kleiman81}, the displayed projectivization is $\Bl_{\YY}(\PP^n\times H)$, proving (ii)(c).

Items (iii)(a) and (iii)(b) follow from \Cref{core_isom}.
For $\mu=\varnothing$ or $\mu=(n^c)$, the construction in that proposition also shows that the universal
kernel $\mathscr K_\mu$ is locally free of rank $q$.
The residual quotient identification and \cite[Proposition 2.2]{Kleiman90} therefore give
\begin{align*}
\QQ_{\bfa}(s,(\mu,1))
&\simeq \Quot^1_{\PP^n\times\QQ_{\bfa}(s,\mu)/\QQ_{\bfa}(s,\mu)}(\mathscr K_\mu)\\
&\simeq \PP_{\PP^n\times\QQ_{\bfa}(s,\mu)}
(\mathscr K_\mu^\vee).
\end{align*}
For $\mu=\varnothing$, the base $\QQ_{\bfa}(s,\mu)$ is $\Gr(s_0,q_0+s_0)$, while for $\mu=(n^c)$, it is $\PP_{\mathbf G}(\mathscr V_c)$.
This proves the remaining case of (iii) and part (iv).
\end{proof}

\begin{remark} When $n=1$ the two trivial cases have simple descriptions. Namely, assume $\lambda=(1^b)$.
If $b=0$, then $\QQ_{\bfa}(s,\varnothing) \simeq\Gr(s_0,q_0+s_0)$ and if $q=1$, then
$\QQ_{\bfa}(r-1,(1^b))\simeq\PP(\mathscr{V}_b)$.
\end{remark}

We end the section with an example of a more general Quot scheme. 
Let $\mathscr{T}_{\PP^2}$ denote the tangent bundle on $\PP^2$.

\begin{thm} \label{connected_tan}
If the Quot scheme $\Quot^P_{\PP^2}(\mathscr{T}_{\PP^2})$ is non-empty, then it is connected.
\end{thm}

\begin{proof}
By twisting, it suffices to consider $\EE=\mathscr{T}_{\PP^2}(-1)$. 
If $\deg(P) = 0$, it is well known that the Quot scheme is irreducible \cite[Theorem~1]{Ellingsrud_Lehn}.
Note that if $\deg(P) = 2$, then every quotient has rank one, whereas if $\deg(P) = 1$, every quotient has a one-dimensional support.

Let $V = S_1$ and consider the Euler sequence
$$
0\longrightarrow\OO_{\PP^2}(-1) \xrightarrow{\sum_i x_i\bfe_i} V^\vee\otimes_\kk\OO_{\PP^2}
\longrightarrow\EE\longrightarrow0
$$
where $\bfe_0,\bfe_1,\bfe_2$ is the basis of $V^\vee$ dual to $x_0,x_1,x_2$. 
We also write $\bfe_i$ for its image in $\EE$.
Given $g \in \GL(3)$, we act on the $x_i$ by left multiplication and on $V^\vee$ by $g^{-T}$. 
Since this preserves $\sum_i x_i\bfe_i$, the action descends to $\EE$ and induces an action on its Quot scheme.

Let the $\mathrm{B}$ be the Borel subgroup of $\GL(3)$ consisting of the upper-triangular matrices. 
Now consider the Borel fixed subsheaves $\mathscr A_1=\OO_{\PP^2}\bfe_2$ and $\mathscr{A}_2=\OO_{\PP^2}\bfe_1\oplus\OO_{\PP^2}\bfe_2$.
Consider a Borel-fixed point $[\EE \twoheadrightarrow \FF]$ in $\Quot_{\PP^2}(\EE)$ and its associated exact sequence
$$
0 \to \KK \to \EE \to \FF \to 0.
$$
We will show that  $\KK \subseteq \mathscr A_1$ if $\FF$ has rank one and $\KK \subseteq \mathscr A_2$ if $\FF$ has dimension one.

First assume that $\FF$ has rank one. Then $\KK$ also has rank one, and the locus where both sheaves are locally free is a dense open subset.
Since $D(x_0)$ is a single $\mathrm{B}$-orbit, this locus contains $D(x_0)$ and in particular, $p_0 = [1:0:0]$.
Note that the fiber of $\EE$ over the point $p_0=[1:0:0]$ is
$\EE|_{p_0}=\langle\bfe_1,\bfe_2\rangle$.
The stabilizer of $p_0$ in the Borel subgroup consists of the matrices
$$
\left\{
\begin{pmatrix}
a&0&0\\
0&b&c\\
0&0&d
\end{pmatrix} \in \GL(3)\right\}
$$
and the action on the fiber $\EE|_{p_0}$ is by the group
$$
\left\{
\begin{pmatrix}
b^{-1}&0\\
-c/(bd)&d^{-1}
\end{pmatrix}
\in \GL(2)
\right\}
$$
Since $\FF$ is locally free at $p_0$, taking fibers gives
$$
0\longrightarrow \KK|_{p_0} \longrightarrow\EE|_{p_0}
\longrightarrow \FF|_{p_0} \longrightarrow 0.
$$
In particular, $\KK|_{p_0}$ is a one-dimensional subspace invariant under the stabilizer of $p_0$.
It follows that  $\KK(p_0) = \langle\bfe_2\rangle = \mathscr A_1(p_0)$.
By equivariance, $\KK|_{D(x_0)}=\mathscr A_1|_{D(x_0)}$ which means that the induced map $\KK\to\EE/\mathscr A_1$ vanishes on $D(x_0)$.
The Euler sequence gives $\EE/\mathscr A_1\cong \II_{[0:0:1]}(1)$, which is torsion-free. 
Hence, the aforementioned map vanishes everywhere proving that $\KK\subseteq\mathscr A_1$.

Now assume that the dimension of $\FF$ is one.
Its support is $\mathrm{B}$-invariant and therefore contains $V(x_0)$.
Similar to the previous case, we will show that $\KK \to \EE/\mathscr{A}_2$ is zero.
Using the Euler sequence we get
$$
0 \longrightarrow \mathscr{A}_2 \longrightarrow \EE \longrightarrow \OO_{V(x_0)} \longrightarrow 0.
$$
We will show that $\KK\to\OO_{V(x_0)}$ is zero.
Consider $p_1=[0:1:0]$ and note that $\EE|_{p_1} = \langle\bfe_0,\bfe_2\rangle$. 
As in the previous case, the only nonzero proper subspace of the fiber invariant under the stabilizer of $p_1$ is $\langle\bfe_2\rangle$. 
Taking fibers this time only gives a right-exact sequence
$$
\KK|_{p_1} \longrightarrow \EE|_{p_1} \longrightarrow \FF|_{p_1} \longrightarrow 0.
$$
Since $p_1 \in \text{Supp}(\FF)$ we have that $\FF|_{p_1} \ne 0$. 
Thus, the image of $\KK|_{p_1} \to \EE|_{p_1}$ is a proper invariant subspace and is contained in $\langle\bfe_2\rangle$. Thus, its image in $(\EE/\mathscr{A}_2)|_{p_1}$ is zero.
By equivariance, this map is zero on fibers throughout the open orbit of $p_1$ in $V(x_0)$. 
Thus, the induced map $\KK|_{V(x_0)} \to  (\EE/\mathscr{A}_2)|_{V(x_0)}$ vanishes on the open orbit of $p_1$.
Since $\EE/\mathscr{A}_2=\OO_{V(x_0)}$ is torsion-free on the integral scheme $V(x_0)$, the map vanishes everywhere; in particular, $\KK \subseteq \mathscr{A}_2$.

Now if $\deg(P) = 2$, then the quotient has rank one and if $\deg(P) = 1$, the quotient is of dimension one. Thus, in either case, the natural morphism
$$
\Quot_{\PP^2}^{\,P-P_{\EE/\mathscr A_i}}(\mathscr A_i) \longrightarrow \Quot^P_{\PP^2}(\EE),
\qquad [\mathscr{A}_i/\KK] \mapsto [\EE/\KK]
$$
contains all the Borel-fixed points.
By \Cref{Pardue_conn}, the domain is connected and thus its image is connected.  
Every connected component of the codomain is $\mathrm{B}$-invariant, since $\mathrm{B}$ is connected, and contains a Borel-fixed point by the Borel fixed-point theorem.
Thus this connected image meets every connected component, i.e., the codomain is connected, as required.
\end{proof}

\section{Radius of Hilbert and Quot schemes}
\label{radius_sec}

In this section we study the radius of Hilbert and Quot schemes and give the first examples of Hilbert schemes with radii larger than 1.
We develop a method to construct components with no common neighbours by using the theory of TNT algebras from the study of the Hilbert scheme of points.
This allows us, rather surprisingly, to show that the Reeves bound is sharp for the Hilbert scheme of curves. 
We also relate how the distance to the lexicographic component of a Quot scheme changes when one slices by a hyperplane.
Using this we show that the radius can grow linearly in the dimension of the subschemes being parameterized.

\begin{notation}
Let $R=\kk[x_1,\ldots,x_e]$ be a standard graded polynomial ring with $e\geq1$.
Since there will be no ambiguity denote the free module by $E_\bfa=\bigoplus_{i=1}^r R(-a_i)\bfe_i$. 
We call a Hilbert function $\mfh$ \textbf{finite} if it has finite support.
The multigraded Quot scheme $\Quot^{\mfh}_{\mathrm{gr}}(E_\bfa)$ is a projective scheme that parameterizes homogeneous quotients of $E_\bfa$ with Hilbert function $\mfh$ \cite[Section 6.2]{Haiman_Sturmfels__multigraded}. The tangent space to a closed point $[E_\bfa/K]$ is $\Hom_R(K,E_\bfa/K)_0$.
\end{notation}

\begin{notation} For a finite Hilbert function $\mfh$, let $H_{\mfh}(u)=\sum_j\mfh(j)u^j$ denote the Hilbert series and let the Hilbert polynomial of the projective cone be
$$
P_{\mfh,c}(t)=\sum_j\mfh(j)\binom{t-j+c}{c}.
$$
\end{notation}

\begin{definition} Let $c \geq 1$ and $S = R[y_0,\dots,y_c]$ with the $y_i$ having degree one. The \textbf{projective cone} of $M = E_\bfa/K$ is 
$C_c(M) := M_S = M \otimes_R S$. Consider the induced morphism
$$
\Quot_{\mathrm{gr}}^{\mfh}(E_\bfa) \to  \Quot^{P_{\mfh,c}}_{\PP^{e+c}}(\widetilde{E_\bfa\otimes_R S}), \qquad [M] \mapsto [\widetilde{C_c(M)}].
$$
If $\HH$ is an irreducible component of $\Quot_{\mathrm{gr}}^{\mfh}(E_{\bfa})$, we define
$$
C_c(\HH) = \overline{\GL(S_1)\cdot\{[\widetilde{C_c(M)}]:[M]\in\HH\}}_{\text{red}},
$$
where $\GL(S_1)$ acts by linear changes of coordinates in $S$ and fixes the basis elements $\bfe_i$.
\end{definition}

\begin{lemma} \label{dimension_cone_comp}
Let $\mfh$ be a finite Hilbert function and let $\HH$ be an irreducible component of $\Quot_{\mathrm{gr}}^{\mfh}(E_{\bfa})$. 
Then $C_c(\HH) \subseteq \Quot_{\PP^{e+c}}(E_\bfa\otimes_R S)$ is irreducible of dimension $\dim\HH+e(c+1)$.
Moreover, every closed point of $C_c(\HH)$ is a linear change of coordinates of $C_c(M)$ for some $[M]\in\HH$.
\end{lemma}
\begin{proof} 
Let $\mathscr{W}$ be the tautological rank $e$ subbundle on $\Gr(e,S_1)$, and consider the associated bundle $Z = \Fr(\mathscr{W}) \times^{\GL(\ell)} \HH$.
By the associated-bundle construction \cite[Section 4.7]{BKS25}, the morphism $Z\to \Gr(e,S_1)$ is Zariski locally trivial with fiber $\HH$. 
In particular, it is proper and irreducible, with $ \dim Z=\dim \HH+e(c+1)$.

Let $\MM$ be the universal quotient on $\Quot^\mfh_{\text{gr}}(E_\bfa)$.
A $T$-point of $\Fr(\mathscr{W}) \times \HH$ consists of a rank $e$ subbundle
$\mathscr{W}_T \subseteq S_1\otimes_\kk \OO_T$, a trivialization $\phi:\OO_T^e \xrightarrow{\sim}\mathscr{W}_T$, and a map $g:T \to \HH$.
Let $R_T=R\otimes_\kk\OO_T$ and note that the trivialization induces a map $R_T\cong\Sym(\mathscr W_T)\hookrightarrow S_T := S \otimes_\kk \OO_T$.
This allows us to define a map from $\Fr(\mathscr{W}) \times \HH \to \Quot(E_\bfa \otimes S)$ by
$(\mathscr{W}_T,\phi,g) \mapsto g^{\star}\MM \otimes_{R_T}S_T$.
Since $S_T$ is locally a polynomial extension of $R_T$, the resulting quotient is flat over $T$ and commutes with base change; in particular, we have a morphism $\Fr(\mathscr{W}) \times \HH \to \Quot(E_\bfa \otimes S)$. 
The morphism is invariant under the diagonal $\GL(e)$-action and hence descends to a morphism $Z \to \Quot(E_\bfa \otimes S)$. 
By construction, its image is closed, irreducible, and equals $C_c(\mathcal H)$. 
It is straightforward to check that this is injective on $\kk$-points giving $ \dim C_c(\mathcal H)=\dim Z =\dim\mathcal H+\ell(c+1)$. 
This completes the proof.
\end{proof}

In general $C_c(\HH)$ need not be a component.
However, the TNT condition gives a sufficient criterion for this to happen.
For a homogeneous quotient $E_\bfa/K$, the coefficientwise derivative $\partial/\partial x_i$ induces an $R$-linear map $\partial_i: K \to E_\bfa/K$ of degree $-1$.

\begin{definition}[{\cite[\S5]{jelisiejew_sivic}}] \label{module_TNT} 
An Artinian module $M = E_\bfa/K$, or the corresponding point on its Quot scheme, has \textbf{trivial negative tangents (TNT)} if the map
$$
\kk^e\longrightarrow\bigoplus_{j<0}\Hom_R(K,M)_j,
\qquad (b_1,\ldots,b_e) \mapsto \sum_{i=1}^e b_i\partial_i
$$
is an isomorphism.
\end{definition}

\begin{thm} \label{module_conet}
Let $\HH\subseteq\Quot^{\mfh}_{\mathrm{gr}}(E_\bfa)$ be an irreducible component of dimension $d$. 
Let $M = E_\bfa/K$ and $\EE = \widetilde{E_\bfa \otimes_R S}$.
Assume that $[M] \in\HH$ has TNT and $\dim_\kk\Hom_R(K,M)_0=d$. 
Then $C_c(M)$ is a smooth point of $\Quot^{P_{\mfh,c}}_{\PP^{e+c}}(\EE)$, and its unique irreducible component is $C_c(\HH)$. 
This component is generically smooth of dimension $d+e(c+1)$, and every one of its points is a projective cone.
\end{thm}

\begin{proof}
Let $C_c(M) = M_S = M \otimes_R S$ and $K_S = K \otimes_R S$.
By \Cref{dimension_cone_comp} it suffices to show that the tangent space dimension at $[\widetilde{M_S}]$ on the Quot scheme is $d+e(c+1)$.
Since $M_S$ has depth at least two, so does $\Hom(K_S,M_S)$, and it follows that 
$$
\Hom_{\PP^{e+c}}(\widetilde{K_S},\widetilde{M_S}) =  H^0(\PP^{e+c},\mathcal{H}om_{\PP^{e+c}}(\widetilde{K_S},\widetilde{M_S})) = \Hom_S(K_S,M_S)_0.
$$
By flat base change we have 
$$
\begin{aligned}
T_{[C_c(M)]}\Quot_{\PP^{e+c}}(\EE)  = \Hom_{\PP^{e+c}}(\widetilde{K_S},\widetilde{M_S}) & =\Hom_S(K_S,M_S)_0\\
&=\bigoplus_{j\geq0} \Hom_R(K,M)_{-j}
                 \otimes_\kk\kk[y_0,\ldots,y_c]_j\\
&=\Hom_R(K,M)_0 \oplus
  \bigl(\langle \partial_1,\ldots,\partial_e \rangle
                  \otimes_\kk \langle y_0,\ldots,y_c\rangle \bigr).
\end{aligned}
$$
In particular, the dimension of the tangent space is $d+ e(c+1)$, as required.
\end{proof}

We now give a cohomological condition for two components to have no common neighbour. 
This is related to the fiber-full stratification, which fixes the cohomology of every twist \cite{CidRuiz_Ramkumar__Fiber_full}. 
The hard part is finding components that satisfy the hypothesis.

\begin{lemma}\label{TNT_separation}
Let $X$ be a projective scheme and $\EE$ a coherent sheaf.
For any integer $m$, the locus
$$
\mathcal{W}_m=\{[\EE/\KK] \in \Quot_X(\EE) :H^1(X,\KK(m))=0\}
$$
is open in $\Quot_X(\EE)$. 
Assume that $\HH_1,\HH_2$ are irreducible components contained in $\mathcal{W}_m$. 
If the function 
\begin{align} \label{dim_0_fun}
[\EE/\KK] \mapsto h^0(X,\KK(m))
\end{align}
takes different values on the two components, then $\HH_1$ and $\HH_2$ are disjoint and have no common neighbours.
\end{lemma}
\begin{proof}
If $\KK_{\Quot_X(\EE)}$ denotes the universal kernel of $\Quot_X(\EE)$, then $\mathcal{W}_m$ is the locus where $H^1(\KK_{\Quot_X(\EE)}(m))$ vanishes which, by upper-semicontinuity, is open. 
For a closed point $[\EE/\KK]  \in \mathcal{W}_m$ we have 
$$
h^0(X,\KK(m))+h^0(X,(\EE/\KK)(m))=h^0(X,\EE(m)).
$$
Both summands are upper semicontinuous and their sum is constant, so each is locally constant. 
In particular, the different values on $\HH_1$ and $\HH_2$ imply that these components are disjoint.

If an irreducible component $\HH$ met both $\HH_1$ and $\HH_2$, then $\HH\cap\mathcal W_m$ would be a nonempty irreducible open subset on which $h^0(X,\KK(m))$ took both values. 
This contradicts the local constancy.
\end{proof}

\begin{definition}\label{TNT_pair}
Let $c\geq1$. 
Two distinct finite Hilbert functions $\mfh,\mfh'$ form a \textbf{$c$-TNT pair} (for $E_{\bfa}$) if they are realized by quotients of $E_{\bfa}$ with TNT and $P_{\mfh,c}=P_{\mfh',c}$. 

When $c=1$, we simply call them a \textbf{TNT pair}. The quantities
$$
d(\mfh)=\sum_j\mfh(j),
\qquad \text{ and }
g(\mfh)=1+\sum_j(j-1)\mfh(j),
$$
are called the \textbf{degree} and \textbf{genus} of the pair.
\end{definition}

\begin{cor}\label{TNT_pair_cor}
Let $(\mfh,\mfh')$ be a $c$-TNT pair for $E_\bfa$ such that their corresponding quotients are smooth on their respective graded Quot schemes. 
Then their cone components in $\Quot^{P_{\mfh,c}}_{\PP^{e+c}}(\EE)$ are disjoint and have no common neighbour. 
In particular,  the radius of $\Quot^{P_{\mfh,c}}_{\PP^{e+c}}(\EE)$ is at least two.
If $c=1$, the Quot scheme has radius exactly two.
\end{cor}
\begin{proof} 
By \Cref{module_conet} a TNT pair for $E_\bfa$ with Hilbert functions $(\mfh,\mfh')$ correspond to smooth points on distinct irreducible components $\HH \subseteq \Quot^{\mfh}_{\mathrm{gr}}(E_\bfa)$ and $\HH' \subseteq \Quot^{\mfh'}_{\mathrm{gr}}(E_\bfa)$, respectively. 
Let $S = R[y_0,\dots,y_c]$ and by \Cref{dimension_cone_comp}, every closed point of  $C_c(\HH)$ or $C_c(\HH')$ is, after a linear change of coordinates, of the form $M_S = (E_\bfa \otimes_R S)/K_S $ or $(E_\bfa \otimes_R S)/K'_S$, respectively.
We have an exact sequence
$$
0 \to K_S \to E_\bfa\otimes_R S \to (E_\bfa\otimes_R S)/K_S \to 0.
$$
Since $\mathrm{depth}_S(K_S) = c+2 \geq 3$, taking the long exact sequence in local cohomology gives $K_S = \bigoplus_{i \in \ZZ}H^0(\widetilde{K_S}(i))$ and $H^1(\widetilde{K_S}(m)) = 0$ for all $m$. 
Similarly, $K'_S =\bigoplus_{i \in \ZZ} H^0(\widetilde{K'_S}(i))$ and $H^1(\widetilde{K'_S}(m)) = 0$ for all $m$. 
Thus, the components $C_c(\HH), C_c(\HH')$ lie in $\mathcal{W}_m$ for all $m$ (see \Cref{TNT_separation}).
The Hilbert series of the quotients on the two components are $\frac{H_{\mfh}(u)}{(1-u)^{c+1}}$ and $\frac{H_{\mfh'}(u)}{(1-u)^{c+1}}$. 
Since $\mfh$ and $\mfh'$ are different, these series are different. 
In particular, the kernels have different sections in some degree $m$.
This implies that \Cref{dim_0_fun} takes different values on $C_c(\HH)$ and $C_c(\HH')$. 
By \Cref{TNT_separation}, the distance between these components is at least three, so the radius of the Quot scheme is at least two. 
If $c=1$, then \Cref{radius_bound} shows that the radius is exactly two, as required.
\end{proof}

\subsection{The two TNT pairs}
To showcase our method, we construct two TNT pairs and use them to give Hilbert schemes
of curves of radius two. 
Although the construction applies more generally to Quot schemes, we restrict to Hilbert schemes for simplicity, since the results are already new in this case.
We use Macaulay's inverse systems to describe
the Artinian algebras. Some standard references are
\cite{Macaulay_inverse_systems,iaCompressed,EliasRossiShortGorenstein,Jel_classifying};
see also \cite[Appendix A]{iakanev} for the divided-power
convention.
\begin{notation}
Let $R=\kk[x_1,\ldots,x_e]$ and let $P=\kk[X_1,\dots,X_e]$. 
We view $P$ as an $R$-module via the \textbf{contraction} action
$$
\bfx^{\boldsymbol{\beta}} \circ \boldsymbol{X}^{\boldsymbol{\alpha}}=
\begin{cases}
\boldsymbol{X}^{\boldsymbol{\alpha}-\boldsymbol{\beta}},&\alpha_i\geq\beta_i\text{ for all }i,\\
0,&\text{otherwise}.
\end{cases}
$$
For homogeneous forms $F_1,\ldots,F_s\in P$, we can consider the annihilator
$\Ann(F_1,\ldots,F_s) = \{f \in R: f \circ F_i=0 \text{ for every }i \}$.
\end{notation}

\begin{lemma} \label{TNT_1} \label{TNT_2} \label{two_TNT_pairs}
The Hilbert functions $(1,6,7)$ and $(1,7,5,1)$ are a TNT pair with degree 14 and genus 7. 
The Hilbert functions $(1,7,10)$ and $(1,8,8,1)$ are also a TNT pair with degree 18 and genus 10.
In all four cases, there is a TNT quotient with this Hilbert function that is smooth on the graded Quot scheme.
\end{lemma}

\begin{proof}
We regard each pair as quotients of the same polynomial ring; extending the ring and quotienting out by the new variable preserves smoothness of the graded Quot schemes and the TNT property.
For $(1,6,7)$, we choose
$$
\begin{aligned}
J=\big(&x_1^2,x_1x_2,x_2^2,x_2x_3,x_3^2,x_3x_6,
        x_4^2,x_5^2,x_5x_6,x_6^2,\\
       &x_1x_3-x_1x_4,\ x_2x_4-x_3x_5,
        \ x_4x_5-x_2x_6,\ x_4x_6-x_1x_6\big)
        \subseteq\kk[x_1,\ldots,x_6].
\end{aligned}
$$
For $(1,7,5,1)$, we choose the ideal
$J'=\Ann(F,Q_1,Q_2)\subseteq\kk[x_1,\ldots,x_7]$,
where
$$
\begin{aligned}
F&=X_1X_2X_3,\\
Q_1&=X_4^2+X_5^2+X_6^2+X_7^2,\\
Q_2&=X_1X_4+X_1X_7+X_2X_5+X_3X_6+X_7^2.
\end{aligned}
$$
For $(1,7,10)$, we choose the ideal 
$$
\begin{aligned}
I=\big(&x_1^2,x_2^2,\ldots,x_7^2,\,
x_1x_2,x_1x_3,x_2x_4,x_2x_6,x_4x_7,x_5x_6,x_5x_7,\\
&x_3x_5-x_4x_6,\,
x_1x_4-x_2x_3,\,
x_1x_6-x_3x_6,\,
x_1x_5-x_6x_7\big) \subseteq \kk[x_1,\ldots,x_7]
\end{aligned}
$$
For $(1,8,8,1)$, we choose the ideal $I'=\Ann(G)$, where
$$
G=X_1X_3X_5+X_2X_4X_6+X_1X_2^2+X_3X_4^2
+X_5X_6^2+X_1X_7X_8+X_3^2X_7+X_5^2X_8
$$
is a divided-power cubic in $\kk[X_1,\dots,X_8]$ (c.f., \cite[Section 4]{Szafarczyk}).

For a homogeneous ideal $\mathfrak{a} \subseteq R$, let $T_j(\mathfrak{a})=\Hom_R(\mathfrak{a},R/\mathfrak{a})_j$. 
A straightforward computation, perhaps using \cite{M2}, gives the following tangent space dimensions:
$$
\begin{array}{c|ccc}
\text{$\mfh$}&\dim T_{-2}&\dim T_{-1}&\dim T_0\\ \hline
(1,6,7)&0&6&98\\
(1,7,5,1)&0&7&67\\
(1,7,10)&0&7&180\\
(1,8,8,1)&0&8&119
\end{array}
$$
Since all the ideals are generated by quadrics, $T_j=0$ for $j<-2$.
In degree $-1$, the derivations are linearly independent and span $T_{-1}$.
Thus, all four algebras have TNT. 

The graded Hilbert scheme with Hilbert function $(1,e,f)$ is the Grassmannian of codimension $f$ subspaces of $R_2$; the corresponding ideal is generated by the subspace along with $R_{\geq 3}$. 
In particular, the maps 
$$
\Gr(7,21) \to \Hilb^{(1,6,7)}(\kk[x_1,\dots,x_6]) \, \text{ and } \, 
\Gr(10,28) \to \Hilb^{(1,7,10)}(\kk[x_1,\dots,x_7])
$$ 
are isomorphisms.
The cubic Gorenstein algebras with Hilbert function $(1,8,8,1)$ form
an open subset of $\PP(\Gamma^3(\kk^8))$, the projective space of divided-power cubics in
eight variables. Since this has dimension $\binom{10}{3}-1$, which is is also the tangent space computed above, we see that $\Hilb^{(1,8,8,1)}(\kk[x_1,\dots,x_8])$ is generically smooth of dimension $119$.

A family of algebras with Hilbert function $(1,7,5,1)$ is obtained as follows.
The choice of a three-dimensional subspace $V \subseteq (\kk^7)^\vee$, a concise divided-power cubic $[F]\in\PP(\Gamma^3(V))$ and a five-dimensional subspace $R_1\circ F \subseteq W \subseteq \Gamma^2((\kk^7)^\vee)$ such that $R_1\circ W=(\kk^7)^\vee$, gives us an algebra $R/\Ann(F,W)$ with Hilbert function $(1,7,5,1)$.
More precisely, let $\mathscr{V}$ be the tautological rank-three subbundle on $\Gr(3,7)$, and consider
the projective bundle $\PP = \PP_{\Gr(3,7)}(\Gamma^3\mathscr{V}).$
Over the open locus of concise cubics, contraction by the universal cubic defines a rank-three subbundle
$\mathscr{D} \subseteq \Gamma^2((\kk^7)^\vee) \otimes \OO_\PP$. Let $\mathscr{Q}$ be the cokernel and note that it has rank $\binom{8}{2}-3 = 25$. 
Finally, these algebras with Hilbert function $(1,7,5,1)$ are parameterized by an open subset of the relative Grassmannian $\Gr_{\PP}(2,\mathscr{Q})$. 
This parameter space is irreducible of dimension $12+9+2(25-2)=67$. 
Since this agrees with the tangent space dimension of $\Hilb^{(1,7,5,1)}(\kk[x_1,\dots,x_7])$, as required.
\end{proof}

We record the two geometric facts that will be used in our construction of the Hilbert scheme.

\begin{definition}
Let $C \subseteq \PP^m$ be an arithmetically Cohen-Macaulay curve with homogeneous coordinate ring $A$.
An \textbf{Artinian reduction} of $A$ is the quotient $A/(\ell_0,\ell_1)$ by a regular sequence of linear forms.
Its Hilbert function is independent of the choice of regular sequence and its sequence of values is called the $h$-\textbf{vector} of $C$.
\end{definition}

\begin{lemma}\label{non_special}
For $(d,g,r)=(14,7,7)$ or $(18,10,8)$, the closure of the locus of smooth genus $g$ curves embedded by a complete nonspecial linear series of degree $d$ in a varying $\PP^r\subseteq\PP^{r+1}$ is a generically smooth irreducible component of the Hilbert scheme.
Its general member is arithmetically Cohen--Macaulay, with $h$-vector $(1,6,7)$ or $(1,7,10)$, respectively.
\end{lemma}

\begin{proof}
Since $d=g+r$, the complete nonspecial embeddings in $\PP^r$ form an open subset of the principal Brill-Noether component \cite[Theorems~19.5 and~19.9]{EH24}. 
Allowing their linear span $\PP^r\subseteq\PP^{r+1}$ to vary therefore gives an irreducible family. 
The image, which consists of nonspecial curves in their complete embedding, corresponds to the locus 
$$
\left\{[C] \in \Hilb^{dt+1-g}(\PP^{r+1}): \text{C is smooth, connected, and } H^1(C,\OO_C(1))= H^1(\PP^{r+1},\II_C(1))=0 \right\}
$$
Indeed, the first vanishing gives $d+1-g=r+1$, while the second vanishing implies $H^0(\PP^{r+1},\OO_{\PP^{r+1}}(1)) \twoheadrightarrow H^0(C,\OO_C(1))$. 
Since both these vanishing conditions are open, the corresponding locus is open in the Hilbert scheme.
By \cite[Theorem 19.5]{EH24}, since $C$ is nonspecial, $H^1(N_{C/\PP^{r+1}})=0$ and we get that the closure of the locus is a generically smooth component. 
By \cite[Corollary 9.2(3)]{eisenbud:syzygies}, the general member is arithmetically Cohen-Macaulay.
Finally, if $C$ is ACM and $\OO_C(1)$ is nonspecial, Riemann-Roch gives $h^0(C,\OO_C(t)) =td+1-g$ for $t \geq 1$. 
Taking second differences of the Hilbert function gives the $h$-vector is $(1,d-g-1,g)$.
Thus the two $h$-vectors are $(1,6,7)$ and $(1,7,10)$, as required.
\end{proof}

\begin{lemma} \label{artinian_reduction}
Let $C\subseteq\PP^m$ be an ACM curve. 
Then $C$ degenerates to the projective cone over any Artinian reduction of its homogeneous coordinate ring.
\end{lemma}
\begin{proof}
We may choose  coordinates $x_1,\ldots,x_{m-1},y_0,y_1$ on  $\PP^m$ so that $y_0,y_1$ are used for the Artinian reduction. 
Let $R = \kk[x_1,\dots,x_{m-1}]$ and $J = (I_C + (y_0,y_1))/(y_0,y_1) \subseteq R$, and let $A = R/J$ be the Artinian reduction.
Consider a torus action obtained by scaling the $x_i$s, namely  $\lambda_t[x_1:\cdots:x_{m-1}:y_0:y_1] = [tx_1:\cdots:tx_{m-1}:y_0:y_1]$.
Let $C_0 = \lim_{t \to 0} \lambda_tC$ be the flat limit of $C$ and note that for every $f\in I_C$, the polynomial $f(x_1,\ldots,x_{m-1},ty_0,ty_1)$ vanishes on $\lambda_tC$.
Thus, at $t=0$, we obtain  $JR[y_0,y_1]\subseteq I_{C_0}$, i.e., $C_0\subseteq\operatorname{Proj}(A[y_0,y_1])$.
Since $y_0,y_1$ form a regular sequence on the coordinate ring, $A[y_0,y_1]$ has the same Hilbert series as $C$. 
By flatness, we must have $C_0 = \Proj(A[y_0,y_1])$, as required.
\end{proof}

We are now ready to give the first example of a Hilbert scheme of curves of radius two. 

\begin{thm}\label{radius_two_14} 
The Hilbert scheme $\Hilb^{14t-6}(\PP^8)$ has radius two. It contains the following generically smooth components:
\begin{enumerate}
\item $\HH_1$, whose general member is a cone over an Artinian algebra with Hilbert function $(1,6,7)$ in a varying $\PP^7$;
\item $\HH_2$, whose general member is a cone over an Artinian algebra with Hilbert function $(1,7,5,1)$;
\item $\mathcal A$, whose general member is a smooth genus-$7$ curve in a varying $\PP^7$;
\item $\mathcal B$, whose general member is a smooth canonical genus-$5$ curve in a $\PP^4$ with two smooth plane cubic tails.
\end{enumerate}
In the last case, the tails are attached at distinct points, the linear span of each cubic meets the $\PP^4$ only at its attaching point, and the three linear spaces span $\PP^8$. The graph induced by these four components is the path
\begin{center}
\begin{tikzpicture}[every node/.style={inner sep=4pt}]
\node (h1) at (0,0) {$\HH_1$};
\node (a) at (2,0) {$\mathcal A$};
\node (b) at (4,0) {$\mathcal B$};
\node (h2) at (6,0) {$\HH_2$};
\draw (h1)--(a)--(b)--(h2);
\end{tikzpicture}
\end{center}
In particular, $d(\HH_1,\HH_2)=3$.
\end{thm}

\begin{proof}
The existence of the  components $\HH_1,\HH_2$ follows from \Cref{two_TNT_pairs,module_conet}.
They are disjoint and have no common neighbour by \Cref{TNT_pair_cor}, which also gives radius two.
The existence of the component $\mathcal A$ follows from \Cref{non_special}.  
A general member of $\mathcal{A}$ is arithmetically Cohen-Macaulay with $h$-vector $(1,6,7)$ \cite[Corollary~8.4]{eisenbud:syzygies}; thus, by \Cref{artinian_reduction}, $\mathcal{A}$ meets $\HH_1$.

We will now show that $\mathcal{B}$ is a component and meets both $\mathcal{A}$ and $\HH_2$.
A general canonical genus $5$ curve is a complete intersection of three quadrics in $\PP^4$. Thus, $\mathcal{B}$ is parameterized by a choice of a $\PP^4$, the three quadrics, the two points on $C_0$, the two $2$-planes, and the cubics through their attaching points. This forms an irreducible family of dimension
$$
\dim\Gr(5,9)+\dim\Gr(3,15)+2 +2\bigl(\dim\Gr(2,8)+8\bigr)=98.
$$
Consider a general point $[C] \in \mathcal{B}$ and let $C=C_0\cup E_1\cup E_2$, with $E_1,E_2$, attached at $p_1,p_2$, respectively. 
Let $N = N_{C/\PP^8}$ be the normal sheaf of $C$ and restricting it to the components gives
\begin{align} \label{C_0res}
0\longrightarrow N_{C_0/\PP^8} \longrightarrow N|_{C_0} \longrightarrow \kk_{p_1}\oplus\kk_{p_2} \longrightarrow0
\end{align}
and
\begin{align} \label{E_res}
0\longrightarrow N_{E_i/\PP^8} \longrightarrow N|_{E_i} \longrightarrow\kk_{p_i} \longrightarrow0,
\end{align}
for $i=1,2$.
We also have
\begin{align} \label{nor_res}
0\longrightarrow N \longrightarrow N|_{C_0}\oplus N|_{E_1}\oplus N|_{E_2} \longrightarrow N|_{p_1}\oplus N|_{p_2} \longrightarrow 0.
\end{align}
Using the normal bundle sequence for the chain $E_i\subseteq\PP^2\subseteq\PP^8$ and twisting it by $-p_i$, we obtain 
$$
0\longrightarrow\OO_{E_i}(3)(-p_i) \longrightarrow N_{E_i/\PP^8}(-p_i) \longrightarrow\OO_{E_i}(1)(-p_i)^{\oplus6} \longrightarrow0.
$$
Since the line bundles on either side have positive degree,  we obtain $H^1(N_{E_i/\PP^8}(-p_i))=0$.
Thus, using \Cref{E_res} we obtain $H^1(N|_{E_i}(-p_i))=0$.
Using the sequence 
$$
0 \longrightarrow N|_{E_i}(-p_i) \longrightarrow N|_{E_i} \longrightarrow N|_{p_i} \longrightarrow 0
$$
we conclude that  $H^0(N|_{E_i})\to N|_{p_i}$ is surjective. Using \Cref{nor_res} we see that 
$H^1(N)\cong H^1(N|_{C_0})$. Taking cohomology of the first restriction sequence now gives
$$
H^0(N|_{C_0})
\longrightarrow\kk^2
\xrightarrow{\delta}
H^1(N_{C_0/\PP^8})
\longrightarrow H^1(N)
\longrightarrow0.
$$
Since $C_0$ is a complete intersection of three quadrics in $\PP^4$, we have $H^1(N_{C_0/\PP^8})\cong H^1(\omega_{C_0})^{\oplus4}\cong\kk^4$.
Under this identification, the connecting homomorphism $\delta$ sends the generators of $\kk_{p_1}$ and $\kk_{p_2}$ to the normal directions of the two elliptic tails modulo $\PP^4$.
These directions are linearly independent since the planes containing the tails, together with $\PP^4$ span $\PP^8$. 
Thus, $\delta$ is injective and $h^1(N)=2$.
By \cite[Chapter XXI, Equation (5.17)]{Arbarello__Geometry_of_algebraic_curves}, for an lci curve $X \subseteq \PP^r$ of degree $d$ and genus $g$, we have $\chi(N_{X/\PP^r})= (r+1)d - (r-3)(g-1)$. 
Using Riemann-Roch, we obtain $h^0(N)=\chi(N)+2 = 96 +2 =98$. 
Thus, $[C] \in \mathcal{B}$ is smooth and $\mathcal{B}$ is an irreducible component.

Consider the exact sequence of homogeneous coordinate rings
$$
0 \to \kk[C] \to \kk[C_0] \oplus (\kk[E_1] \oplus \kk[E_2]) \to \kk[p_1] \oplus \kk[p_2] \to 0. 
$$ 
Since $E_1,E_2,C_0$ are ACM and the rightmost term has depth one, the depth lemma shows that $C$ is ACM. 
Furthermore, we see that the Hilbert series of $\kk[C]$ is
$$
\frac{1+3u+3u^2+u^3+2(1+u+u^2)-2(1-u)}{(1-u)^2}
=\frac{1+7u+5u^2+u^3}{(1-u)^2}.
$$
Choose two linear forms forming a regular sequence on both $C$ and $C_0$, and let $A$ and $A_0$ denote the corresponding Artinian reductions.
The inclusion $C_0\subseteq C$ induces a surjection $A\twoheadrightarrow A_0$, and hence an inclusion of inverse systems $N'\subset N$.
Since $A_0$ is Gorenstein, its inverse system is generated by a concise ternary cubic $F$. 
Moreover, we have $(N')_3=N_3=\langle F\rangle$, $(N')_2=S_1\circ F\subset N_2$ and $\dim_\kk N_2=5$.
Thus, the inverse system of $A$ has cubic part generated by $F$ and a five-dimensional quadratic part containing its contractions. 
It follows that $A$ lies in the closure of the family with Hilbert function $(1,7,5,1)$ constructed in \Cref{two_TNT_pairs}.
Thus, \Cref{artinian_reduction} implies that $\mathcal{B}$ meets $\HH_2$.

It remains to show that $\mathcal A\cap\mathcal B\neq\varnothing$. 
We will prove that both components contain an isomorphic projection of $C$ into $\PP^7$. 
Let $L=\OO_C(1)$ and note that $L|_{C_0}\cong\omega_{C_0}$.
Since each $E_i$ meets $C_0$ transversely at $p_i$, the ideal sheaf sequence for $C_0\subset C$, tensored with $L$, gives
$$
0\longrightarrow L|_{E_1}(-p_1)\oplus L|_{E_2}(-p_2) \longrightarrow L \longrightarrow \omega_{C_0} \longrightarrow 0.
$$
Taking cohomology we have 
$$
0\longrightarrow \kk^4\longrightarrow H^0(C,L) \xrightarrow{\rho} H^0(C_0,\omega_{C_0})\longrightarrow 0
$$ 
since $H^0(L|_{E_i}(-p_i)) = H^0(\OO_{E_i}(2)) = \kk^2$ by Riemann-Roch.

Since the secant variety of $C_0$ is three-dimensional, we can choose a point $q \in \mathrm{span}(C_0) = \PP^4$ but outside the secant variety. 
Let $V\subseteq H^0(C_0,\omega_{C_0})$ be the hyperplane corresponding to $q$ and let $W=\rho^{-1}(V)$.
Let $C'$ denote the projection of $C$ from $q$ into a $\PP^7$. Since $C'$ is isomorphic to $C$, the linear projection is a flat family and we have $[C'] \in \mathcal{B}$.

By Serre duality, $H^1(C_0,\OO_{C_0})\cong H^0(C_0,\omega_{C_0})^\vee$.
We can choose $\eta\in H^1(C_0,\OO_{C_0})$ corresponding to a nonzero linear functional with kernel $V$.
Since $\mathrm{Pic}^8(C_0)$ is smooth with tangent space $H^1(C_0,\OO_{C_0})$, we can find a deformation of $\omega_{C_0}$ over a smooth curve $T$ lifting the tangent vector $\eta$. Since $\omega_{C_0}$ is the unique special line bundle of degree eight on $C_0$, the nonzero tangent direction $\eta$
ensures that one-parameter family is nonspecial for $t\ne0$.
By leaving the restrictions of $\OO_C(1)$ to the elliptic tails fixed, we obtain a non-trivial one-parameter family of line bundles on $C \times T$ of degree fourteen with $\mathscr{L}_0 = \OO_C(1)$.

For $t\ne 0$, we have $H^1(C,\mathscr{L}_t)=0$ and $h^0(C,\LLL_t)=8$. 
Let $\pi:C\times T \to T$ be the projection and since $\pi_{\star} \LLL$ is locally free of rank eight, we have an injective map $(\pi_{\star}\LLL)\otimes \kk(0) \to H^0(C,\LLL)$, where $\kk(0)$ is the residue field of $0 \in T$.
Its image, denoted by $W'$, is eight dimensional. By \cite[Proposition 3.3.4]{Sernesi__Deformations}, the obstruction to extending a section of $\OO_C(1)$ to first order
is the cup-product map
$
H^0(C,\OO_C(1))\xrightarrow{\rho}H^0(C_0,\omega_{C_0})
\xrightarrow{\smile\eta}H^1(C_0,\omega_{C_0}).
$
Since every section in $W'$ of course extends to first order, we have  $W'\subseteq\rho^{-1}(V)=W$. Since both spaces have dimension eight, we have $W'=W$.

Choosing a basis of sections of $\pi_\star \mathscr{L}$ over a neighbourhood $T' \ni 0$, we obtain a family of morphisms $C\times T' \to\PP^7$ whose special fiber is the embedding defined by $W$, namely $C'$.
Since the special fiber is a closed immersion, we can shrink $T'$ to obtain a family of closed immersions.
Thus, $C'$ is a flat limit of nonspecially embedded nodal curves in $\PP^7$.
Restricting the Euler sequence of $\PP^7$ to $C_t$ gives us
$$
0 \to \OO_{C_t} \to \OO_{C_t}(1)^{\oplus 8} \simeq \LLL_t^{\oplus 8} \to \mathscr{T}_{\PP^7}|_{C_t} \to 0.
$$
Since $\LLL_t$ is nonspecial, we can take cohomology and obtain $H^1(C_t,\mathscr{T}_{\PP^7}|_{C_t})=0$.
Since $C$ is a nodal curve, \cite[Lemma 2.2]{HRS04} now implies that the nodes can be smoothed
inside $\PP^7$.
The resulting smooth curve must have degree fourteen and genus seven.
By upper semicontinuity, the smoothing remains nonspecial, and hence belongs to $\mathcal{A}$ by \Cref{non_special}.
In particular, $[C_t] \in \mathcal{A}$ for $t \ne 0$, which in turn implies that $[C'] \in \mathcal{A}$, as required.
\end{proof}

In a similar, but easier, fashion we can give another example of a Hilbert scheme of curves of radius two.

\begin{thm}\label{radius_two_18}
The Hilbert scheme $\Hilb^{18t-9}(\PP^9)$ has radius two.
It contains the following generically smooth components:
\begin{enumerate}
\item $\HH_1$, whose general member is a cone over an Artinian
algebra with Hilbert function $(1,7,10)$ in a varying $\PP^8$;
\item $\HH_2$, whose general member is a cone over a Gorenstein
algebra having Hilbert function $(1,8,8,1)$;
\item $\mathcal A$, whose general member is a smooth genus-$10$
curve in a varying $\PP^8$;
\item $\mathcal B$, whose general member is a smooth canonical
genus-$10$ curve in $\PP^9$.
\end{enumerate}
The graph induced by these four components is the path
\begin{center}
\begin{tikzpicture}[every node/.style={inner sep=4pt}]
\node (h1) at (0,0) {$\HH_1$};
\node (a) at (2,0) {$\mathcal A$};
\node (b) at (4,0) {$\mathcal B$};
\node (h2) at (6,0) {$\HH_2$};
\draw (h1)--(a)--(b)--(h2);
\end{tikzpicture}
\end{center}
In particular, $d(\HH_1,\HH_2)=3$.
\end{thm}

\begin{proof}
The existence of the  components $\HH_1,\HH_2$ follows from \Cref{two_TNT_pairs,module_conet}.
They are disjoint and have no common neighbour by \Cref{TNT_pair_cor}, which also gives radius two.
The existence of the component $\mathcal A$ follows from \Cref{non_special}. 
The existence of $\mathcal{B}$ follows from the fact that the locus of smooth canonically embedded curves of genus $g$ is open in $\Hilb^{(2g-2)t+1-g}(\PP^{g-1})$. 
Indeed, by Riemann-Roch and Serre duality, every smooth nondegenerate curve with this Hilbert polynomial is canonically embedded. It is irreducible and generically smooth by \cite[Example 3.2.5(vi)]{Sernesi__Deformations}. 

The general member of $\mathcal{A}$ is arithmetically Cohen-Macaulay with $h$-vector $(1,7,10)$ \cite[Corollary 8.4]{eisenbud:syzygies}; thus, by \Cref{artinian_reduction}, $\mathcal{A}$ meets $\HH_1$. By \cite[Section 9.1]{eisenbud:syzygies}, a general canonical curve of genus 10 is arithemetically Gorenstein with $h$-vector $(1,8,8,1)$; thus, by \Cref{artinian_reduction}, $\mathcal{B}$ meets $\HH_1$. 
Finally, $\mathcal{A}$ meets $\mathcal{B}$ by repeating the same deformation of line bundle argument given in \Cref{radius_two_14}.
\end{proof}

\subsection{Unbounded radii}
We use hyperplane sections to compare distances between components of Quot schemes. 
We show that taking a cone increases every positive lex-distance, and consequently, gives Quot schemes with arbitrarily large radius. 
Since the rest of the paper will involve going back and forth between two Quot schemes, it will be simpler to fix some notation for the rest of the paper.

\begin{notation}
Let $n\geq2$, $\QQ=\QQ_{\bfa}(s,\lambda)$ and $P=P_{\bfa,s,\lambda}$.
Let $\LL\subseteq\QQ$ denote the lexicographic component and $\ell$ the lexicographic point (see \Cref{lex_smooth}).
Fix the hyperplane $H=V(x_n)$ and let 
$$
\QQ_H=
\Quot_H^{\Delta P}(\EE_{\bfa}|_H)
=\QQ_{\bfa,H}(s,\partial\lambda),
$$
where $\Delta P(t)=P(t)-P(t-1)$ and $H$ denotes the ambient projective space. Let $\ell_H$ denote the lexicographic point in $\QQ_H$.
\end{notation}

\begin{definition}\label{def:hyperplane-regular-locus}
The \textbf{hyperplane-regular locus} is
$$
\UU_H=
\left\{ [\EE_{\bfa}\twoheadrightarrow\FF]\in\QQ \ \middle|\ x_n:\FF(-1) \to \FF \text{ is injective} \right\}.
$$
\end{definition}

\begin{lemma}\label{lem:hyperplane-slicing}
The locus $\UU_H$ is open and contains every Borel-fixed point.
Restriction to $H$ induces a morphism
$$
\rho_H:\UU_H\longrightarrow\QQ_H,
\qquad [\EE_{\bfa}\twoheadrightarrow\FF] \mapsto [\EE_{\bfa}|_H\twoheadrightarrow\FF|_H].
$$
Moreover, if $z\in\LL$ is Borel-fixed, then $\rho_H(z)=\ell_H$.
\end{lemma}

\begin{proof}
Let $\pi:\PP^n\times\QQ\to\QQ$ be the projection, let $\FF_\lambda$ be the universal quotient, and consider
$\mu:\FF_\lambda(-1) \xrightarrow{x_n} \FF_\lambda$.
By the fiberwise injectivity criterion \cite[\href{https://stacks.math.columbia.edu/tag/00ME}{Tag 00ME}]{stacks_project}, $\mu$ is injective on the fiber at a point of $\PP^n\times\QQ$ if and only if its stalk is injective and its cokernel is flat over $\QQ$
there. 
These two conditions are both open and thus define an open subset $W \subseteq \PP^n \times \QQ$. Indeed, the first is the complement of $\Supp(\ker\mu)$, and the second is open by \cite[\href{https://stacks.math.columbia.edu/tag/0399}{Tag 0399}]{stacks_project}.
Since $\pi$ is proper, the locus $ \UU_H=\QQ\setminus \pi((\PP^n\times\QQ)\setminus W)$ is open.

Over this open subset, the same criterion gives an exact sequence
$$
0 \longrightarrow \FF_\lambda(-1) \xrightarrow{x_n} \FF_\lambda \longrightarrow \FF_\lambda|_{H\times\UU_H} \longrightarrow 0
$$
whose cokernel is flat over $\UU_H$.
The fibers of the cokernel have Hilbert polynomial $\Delta P = P(t) - P(t-1)$.
In particular, the morphism $\rho_H$ is precisely obtained by restricting the universal quotient to $H \times \UU_H$. By \Cref{Pardue_Borel}, $x_n$ is a nonzero divisor on $E_{\bfa}/K$ for every saturated Borel-fixed submodule $K$. In particular, $\UU_H$ contains all the Borel-fixed points.

Let $K \subseteq E_\bfa$ be a Borel-fixed submodule on the lexicographic component. 
Let $S'=\kk[x_0,\ldots,x_{n-1}]$ and $E_{S'}=\bigoplus_iS'(-a_i)\bfe_i$. 
Since the minimal generators of $K$ do not involve $x_n$, the saturated kernel of its
slice is
$$
(K\cap E_{S'})^{\mathrm{sat}} =(K\cap E_{S'}):x_{n-1}^{\infty} =(K:x_{n-1}^{\infty})\cap E_{S'}.
$$
By \Cref{double_saturation} we have $K:x_{n-1}^{\infty} = \ell$. It follows that $(K\cap E_{S'})^{\mathrm{sat}} = (\ell \cap E_{S'})^{\mathrm{sat}} = \ell_H$, implying that $\rho_H([E_{\bfa}/K])=\ell_H$, as required.
\end{proof}

\begin{definition}\label{not:component-slice}
For an irreducible component $\HH\subseteq\QQ$, define its \textbf{slice image} by
$$
\rho_H(\HH) := \overline{\rho_H(\HH\cap\UU_H)}\subseteq\QQ_H.
$$
\end{definition}

\begin{remark}
Since every irreducible component $\HH \subseteq \QQ$ contains a Borel-fixed point, it follows that $\rho_H(\HH)$ is always a non-empty irreducible closed subset. This is usually not a component.
\end{remark}

\begin{thm}[Slicing criterion]\label{slicing_criterion}
Let $z\in\UU_H$ and suppose $d_{\lex}(\rho_H(z))>0$. 
Then
$$
d_{\lex}(z) \geq 1+d_{\lex}(\rho_H(z)).
$$
\end{thm}

\begin{proof}
Let $k=d_{\lex}(z)$, and choose a chain of components $\HH_0,\HH_1,\ldots,\HH_k$ attaining this distance, where $\HH_0$ is the lexicographic component and $z \in \HH_k$.
If $k=0$, the irreducible closed subset $\rho_H(\LL)$ contains both $\rho_H(z)$ and the lexicographic point $\ell_H$. These points therefore lie on a common component, contradicting the hypothesis.

Assume $k \geq 1$ and for each $1\leq i\leq k$, choose an irreducible component $\mathcal{X}_i$ of $\QQ_H$ containing $\rho_H(\HH_i)$. 
To see that $\XX_1,\dots,\XX_k$ is a chain of components that meet, note that the intersection $\HH_{i-1} \cap \HH_i$ always contains Borel-fixed points. 
By \Cref{lem:hyperplane-slicing}, these Borel-fixed points lie in $\UU_H$, which in turn implies that they lie in $\XX_{i-1} \cap \XX_i$. 
By \Cref{lem:hyperplane-slicing} we have $\rho_H(z_1)=\ell_H \in \XX_1$. 
In particular, $\XX_1,\ldots,\XX_k$ gives a walk of length at most $k-1$ between components containing $\ell_H$ and $\rho_H(z)$. 
In particular, $d_{\lex}(\rho_H(z))\leq k-1$, as required.
\end{proof}

\begin{notation}
Consider a closed point $z = [E_{\bfa}/K]$ on $\Quot_{\PP^n}(\EE_\bfa)$. 
For $m \geq 0$, let $S^{[m]}= S[x_{n+1},\dots,x_{n+m}]$ and let $C_m(z)$ denote the quotient obtained by sheafifiying $(E_\bfa/K) \otimes_S S^{[m]}$. 
Thus $C_0(z) =z$, and for $m\geq 1$ this is the projective cone with vertex $\PP^{m-1}$. 
Let $P_m$ be its Hilbert polynomial and let
$$
\QQ_m=\Quot_{\PP^{n+m}}^{P_m}
\left(\bigoplus_i\OO_{\PP^{n+m}}(-a_i)\right).
$$
For $m\geq 1$, $x_{n+m}$ is a non-zero divisor on $C_m(z)$ and thus, the corresponding hyperplane section is $C_m(z)$. 
In particular, $\Delta P_m=P_{m-1}$.
\end{notation}

\begin{cor}\label{cone_distance}
With notation as above, if $d=d_{\lex}(z)>0$ then $d_{\lex}(C_m(z))\geq d+m$  and, in particular,
$$
\rad(\QQ_m)\geq
\left\lceil\frac{d+m}{2}\right\rceil.
$$
\end{cor}

\begin{proof}
For $m\geq1$, multiplication by $x_{m+n}$ is injective on the cone quotient, and restriction to $V(x_{m+n})$ gives $C_{m-1}(z)$. Applying \Cref{slicing_criterion} successively gives
$$
d_{\lex}(C_m(z)) \geq1+d_{\lex}(C_{m-1}(z)) \geq d+m.
$$

Every component containing $C_m(z)$ therefore has distance at least $d+m$ from every component containing the lexicographic point. 
Thus the incidence graph has diameter at least $d+m$. 
Since the diameter of a connected graph is at most twice its radius, the result follows.
\end{proof}

This corollary allows us to easily generate a number of examples with large radii. Here is one starting with a non-smoothable algebra in $\PP^4$.

\begin{cor} Let $
P_m(t) = \binom{t+m}{m} +4\binom{t+m-1}{m} +3\binom{t+m-2}{m}.
$
Then the Hilbert scheme $\Hilb^{P_m}(\PP^{m+4})$ contains cones over algebras of length $8$ and we have $$
\rad\bigl(\Hilb^{P_m}(\PP^{m+4})\bigr) \geq\left\lceil\frac{m+1}{2}\right\rceil.
$$
\end{cor}
\begin{proof} 
Let $S =\kk[x_0,\dots,x_4]$ and consider $Z \subseteq \PP^4$ defined by
$$
J=(x_1^2,x_1x_2,x_2^2,x_3^2,x_3x_4,x_4^2, x_1x_4+x_2x_3) \subseteq S.
$$
The algebra $A = \kk[x_1,\dots,x_4]/(J+x_0)$ has Hilbert function $(1,4,3)$ and is not smoothable \cite[Exc. 5.8]{HarDeform}. 
It follows that $[Z]$ does not lie on the lexicographic component, i.e., $d_{\lex}([Z]) = 1$.
The result now follows from \Cref{cone_distance}.
\end{proof}

\bibliographystyle{amsalpha} 
\bibliography{references}

\end{document}